\documentclass{elsarticle}
\usepackage{hyperref}
\usepackage[utf8]{inputenc}
\usepackage{wrapfig}
\usepackage{amsmath}
\usepackage{amsthm}
\usepackage{amssymb}
\usepackage[tmargin=1in,bmargin=1in]{geometry}
\usepackage{caption}
\usepackage{graphicx}
\usepackage{latexsym}
\usepackage{pdfsync}
\usepackage[boxed]{algorithm}
\usepackage{algorithm}
\usepackage{algpseudocode}
\usepackage{multirow}
\usepackage{rotating}
\usepackage{color}
\usepackage{url}
\usepackage{subfigure}
\usepackage{bigstrut}
\usepackage{bigints}
\usepackage[normalem]{ulem}

\usepackage[normalem]{ulem}
\usepackage{hyperref}

\usepackage{booktabs}
\usepackage{siunitx}

\hypersetup{
    colorlinks=true,              
    linkcolor=blue,               
    citecolor=blue,               
    filecolor=black,              
    urlcolor=red,               
    bookmarks=false,
    pdffitwindow=true,
    pdfpagelayout=SinglePage
}

\definecolor{mypink}{cmyk}{0, 0.7808, 0.4429, 0.1412}

\definecolor{mybrown}{cmyk}{0, 0.60, 0.95, 0.63}
\definecolor{darkgreen}{cmyk}{0.98, 0, 0.36, 0.22}

\NewDocumentCommand{\tens}{t_}
 {%
  \IfBooleanTF{#1}
   {\tensop}
   {\otimes}%
 }

\newcommand{\diplomacy}[1]{\iffalse #1 \fi}
\newcommand{\etal}{\textit{et al.\ }}

\newcommand{\be}{\begin{equation}}
\newcommand{\ee}{\end{equation}}
\newcommand{\ben}{\begin{equation*}}
\newcommand{\een}{\end{equation*}}
\newcommand{\bea}{\begin{eqnarray}}
\newcommand{\eea}{\end{eqnarray}}
\newcommand{\bean}{\begin{eqnarray*}}
\newcommand{\eean}{\end{eqnarray*}}

\begin{document}

\title{Solving Inverse-PDE Problems with Physics-Aware Neural Networks}

\author[MECHE]{Samira Pakravan\corref{cor1}\fnref{fn1}}
\author[MECHE]{Pouria A. Mistani\fnref{fn2}}
\author[ASTRO]{Miguel A. Aragon-Calvo}
\author[MECHE,CS]{Frederic Gibou}
\address[MECHE]{Department of Mechanical Engineering, University of California, Santa Barbara, CA 93106-5070}
\address[CS]{Department of Computer Science, University of California, Santa Barbara, CA 93106-5110}
\address[ASTRO]{Instituto de Astronom\'ia, UNAM, Apdo. Postal 106, Ensenada 22800, B.C., M\'exico}
\cortext[cor1]{Corresponding author: spakravan@ucsb.edu}
\fntext[fn1,fn2]{Equal contribution.}

\begin{abstract}
We propose a novel composite framework to find unknown fields in the context of inverse problems for partial differential equations (PDEs). We blend the high expressibility of deep neural networks as universal function estimators with the accuracy and reliability of existing numerical algorithms for partial differential equations as custom layers in semantic autoencoders. Our design brings together techniques of computational mathematics, machine learning and pattern recognition under one umbrella to incorporate domain-specific knowledge and physical constraints to discover the underlying hidden fields. The network is explicitly aware of the governing physics through a hard-coded PDE solver layer in contrast to most existing methods that incorporate the governing equations in the loss function or rely on trainable convolutional layers to discover proper discretizations from data. This subsequently focuses the computational load to only the discovery of the hidden fields and therefore is more data efficient. We call this architecture Blended inverse-PDE networks (hereby dubbed BiPDE networks) and demonstrate its applicability for recovering the variable diffusion coefficient in Poisson problems in one and two spatial dimensions, as well as the diffusion coefficient in the time-dependent and nonlinear Burgers' equation in one dimension. We also show that this approach is robust to noise.
\end{abstract}


\maketitle

\section{Introduction}
\label{sec::introduction}
Inverse differential problems, where given a set of measurements one seeks a set of optimal parameters in a governing differential equation, arise in numerous scientific and technological domains. Some well-known applications include X-ray tomography \cite{epstein2007introduction,natterer2001mathematics}, ultrasound \cite{van2019deep}, MRI imaging \cite{jin2017deep}, and transport in porous media \cite{J.-P.-Fouque:2007aa}. Moreover, modeling and control of dynamic complex systems is a common problem in a broad range of scientific and engineering domains, with examples ranging from understanding the motion of bacteria colonies in low Reynolds number flows \cite{samsami2020stability}, to the control of spinning rotorcrafts in high speed flights \cite{hedayatpour2017unified,hedayatpour2019precision}. Other applications in medicine, navigation, manufacturing, \textit{etc.} need estimation of the unknown parameters in \textit{real-time}, \textit{e.g.} in electroporation \cite{zupanic2012treatment,mistani2019parallel} the pulse optimizer has to be informed about tissue parameters in microsecond time. On the other hand, high resolution data-sets describing spatiotemporal evolution of complex systems are becoming increasingly available by advanced multi-scale numerical simulations (see \emph{e.g.} \cite{mistani2019parallel, mistani2018island}). These advances have become possible partly due to recent developments in discretization techniques for nonlinear partial differential equations with sharp boundaries (see \emph{e.g.} the reviews \cite{Gibou:2019aa, Gibou:2018aa}). However, solving these inverse problems poses substantial computational and mathematical challenges that makes it difficult to infer reliable parameters from limited data and in real-time.

The problem can be mathematically formulated as follows. Let the values of $u=u(t, x_1, \ldots, x_n)$ be given by a set of measurements, which may include noise. Knowing that $u$ satisfies the partial differential equation:
\begin{eqnarray*}
\frac{\partial u}{\partial t} = f \left (t, x_1, \ldots, x_n; u, \frac{\partial u}{\partial x_1}, \ldots \frac{\partial u}{\partial x_n}; \frac{\partial^2 u}{\partial x_1 \partial x_1}, \ldots \frac{\partial^2 u}{\partial x_1 \partial x_n}; \ldots; \mathbf{c} \right),
\end{eqnarray*}
find the hidden fields stored in $\mathbf{c}$, where the hidden fields can be constant or variable coefficients (scalars, vectors or tensors).

Deep neural networks have, rather recently, attracted considerable attention for data modeling in a vast range of scientific domains, in part due to  freely available modern deep learning libraries (in particular \texttt{TensorFlow} \cite{abadi2016tensorflow}). For example, deep neural networks have shown astonishing success in emulating sophisticated simulations \cite{he2019learning,zhang2019dark,zamudio2019higan,chandrasekaran2019solving,sinitskiy2018deep}, discovering governing differential equations from data \cite{Raissi:2017aa,berg2019data,long2018pde,schaeffer2017learning}, as well as potential applications to study and improve simulations of multiphase flows \cite{Gibou:2019aa}. We refer the reader to \cite{owhadi2019operator,owhadi2019statistical} for a comprehensive survey of interplays between numerical approximation, statistical inference and learning. However, these architectures require massive datasets and extensive computations to train numerous hidden weights and biases. Therefore, reducing complexity of deep neural network architectures for inverse problems poses a significant practical challenge for many applications in physical sciences, especially when the collection of large datasets is a prohibitive task \cite{Raissi:2018aa}. One remedy to reduce the network size is to embed the knowledge from existing mathematical models \cite{Stinis:2019aa} or known physical laws within a neural network architecture \cite{Ling:2016aa,Geng:2019aa}. Along these lines, semantic autoencoders were recently proposed by Aragon-Calvo \cite{aragon2019self}, where they replaced the decoder stage of an autoencoder architecture with a given physical law that can reproduce the provided input data given a physically meaningful set of parameters. The encoder is then constrained to discover optimal values for these parameters, which can be extracted from the bottleneck of the network after training. We shall emphasize that this approach reduces the size of the unknown model parameters, and that the encoder can be used independently to infer hidden parameters in real time, while adding interpretability to deep learning frameworks. Inspired by their work, we propose to blend traditional numerical solver algorithms with custom deep neural network architectures to solve inverse PDE problems more efficiently, and with higher accuracy.

\subsection{Existing works}
Recently, the most widely used approach for solving forward and inverse partial differential equations using neural networks has been the constrained optimization technique. These algorithms augment the cost function with terms that describe the PDE, its boundary and its initial conditions, while the neural network acts as a surrogate for the solution field. Depending on how the derivatives in the PDEs are computed, there may be two general classes of methods that we review in the next paragraph.

In the first class, spatial differentiations in the PDE are performed exclusively using automatic differentiation, while temporal differentiation may be handled using the traditional Runge-Kutta schemes (called \textit{discrete time models}) or using automatic differentiations (called \textit{continuous time models}) \cite{raissi2017physics}. In these methods, automatic differentiation computes gradients of the output of a neural network with respect to its input variables. Hence, the input must always be the independent variables, \emph{i.e.} the input coordinates $\mathbf{x}$, time and the free parameters. In this regard, network optimization aims to calibrate the weights and biases such that the neural network outputs the closest approximation of the solution of a PDE; this is enforced through a regularized loss function. An old idea that was first proposed by Lagaris \etal (1998) \cite{lagaris1998artificial}. In 2015, the general framework of solving differential equations as a learning problem was proposed by Owhadi \cite{owhadi2015bayesian,Owhadi2015,owhadi2017multigrid} which revived interest in using neural networks for solving differential equations in recent years. Raissi \etal (2017) \cite{raissi2017physics,Raissi2017PhysicsID} presented the regularized loss function framework under the name \textit{physics informed neural networks} or PINNs and applied it to time-dependent PDEs. Ever since, other authors have mostly adopted PINNs, see \emph{e.g.} \cite{Sirignano:2018aa,bar2019unsupervised}. The second class of constrained optimization methods was proposed by Xu and Darve \cite{xu2019neural} who examined the possibility of directly using pre-existing finite discretization schemes within the loss function.

An alternative approach for solving PDE systems is through explicit embedding of the governing equations inside the architecture of deep neural networks via convolutional layers, activation functions or augmented neural networks. Below we review some of these methods:
\begin{itemize}
\item A famous approach is PDE-Net \cite{long2018pde,long2019pde} which relies on the idea of numerical approximation of differential operators by convolutions. Therefore, PDE-Nets use convolution layers with trainable and constrained kernels that mimic differential operators (such as $\rm U_x, U_y, U_{xx}, \cdots$) whose outputs are fed to a (symbolic) multilayer neural network that models the nonlinear response function in the PDE system, \textit{i.e.} the right hand side in $\rm U_t = F(U, U_x, U_y, U_{xx}, \cdots)$. Importantly, PDE-Nets can only support \textit{explicit} time integration methods, such as the forward Euler method \cite{long2018pde}. Moreover, because the differential operators are being learned from data samples, these methods have hundreds of thousands of trainable parameters that demand hundreds of data samples; \textit{e.g.} see section 3.1 in \cite{long2018pde} that uses $20$ $\delta t$-blocks with $17,000$ parameters in each block, and use $560$ data samples for training. 

\item Berg and Nystr{\"o}m \cite{berg2017neural} (hereby BN17) proposed an augmented design by using neural networks to estimate PDE parameters whose output is fed into a forward finite element PDE solver, while the adjoint PDE problem is employed to compute gradients of the loss function with respect to weights and biases of the network using automatic differentiation. Even though their loss function is a simple $\rm L_2$-norm functional, the physics is not localized in the structure of the neural network as the adjoint PDE problem is also employed for the optimization process. It is important to recognize that in their approach the numerical solver is a separate computational object than the neural network, therefore computing gradients of error functional with respect to the network parameters has to be done explicitly through the adjoint PDE problem. Moreover, their design can not naturally handle trainable parameters in the numerical discretization itself, a feature that is useful for some meshless numerical schemes. In contrast, \emph{in BiPDEs the numerical solver is a computational layer added in the neural network architecture and naturally supports trainable parameters in the numerical scheme.} For example in the meshless method developed in section \ref{sec::meshfree} we leverage this unique feature of BiPDEs to also train for shape parameters and interpolation seed locations of the numerical scheme besides the unknown diffusion coefficient. 

\item Dal Santos \etal \cite{dal2020data} proposed an embedding of a reduced basis solver as \textit{activation function} in the last layer of a neural network. Their architecture resembles an autoencoder in which the decoder is the reduced basis solver and the parameters at the bottleneck ``are the values of the physical parameters themselves or the affine decomposition coefficients of the differential operators'' \cite{dal2020data}. 

\item Lu \etal \cite{lu2020extracting} proposed an unsupervised learning technique using variational autoencoders to extract physical parameters (not inhomogeneous spatial fields) from noisy spatiotemporal data. Again the encoder extracts physical parameters and the decoder propagates an initial condition forward in time given the extracted parameters. These authors use convolutional layers both in the encoder to extract features as well as in the decoder with recurrent loops to propagate solutions in time; \textit{i.e.} the decoder leverages the idea of estimating differential operators with convolutions. Similar to PDE-Nets, this architecture is also a ``PDE-integrator with explicit time stepping'', and also they need as few as 10 samples in the case of Kuramoto-Sivashinsky problem.
\end{itemize} 
In these methods, a recurring idea is treating latent space variables of autoencoders as physical parameters passed to a physical model decoder. This basic idea pre-dates the literature on solving PDE problems and has been used in many different domains. Examples include Aragon-Calvo \cite{aragon2019self} who developed a galaxy model fitting algorithm using \textit{semantic autoencoders}, or Google Tensorflow Graphics \cite{tensorflowGraphics} which is a well-known application of this idea for scene reconstruction.

\subsection{Present work}
Basic criteria of developing numerical schemes for solving partial differential equations are \textit{consistency} and \textit{convergence} of the method, \textit{i.e.} increasing resolution of data should yield better results. Not only there is no guarantee that approximating differential operators through learning convolution kernels or performing automatic differentiations provide a consistent or even stable numerical method, but also the learning of convolution kernels to approximate differential operators requires more data and therefore yield less data-efficient methods. Therefore it seems reasonable to explore the idea of blending classic numerical discretization methods in neural network architectures, hence informing the neural network about proper discretization methods. This is the focus of the present manuscript.

In the present work, we discard the framework of constrained optimization altogether and instead choose to explicitly blend fully traditional finite discretization schemes as the decoder layer in semantic autoencoder architectures. In our approach, the loss function is only composed of the difference between the actual data and the predictions of the solver layer, but contrary to BN17 \cite{berg2017neural} we do not consider the adjoint PDE problem to compute gradients of the error functional with respect to network parameters. This is due to the fact that in our design the numerical solver is a custom layer inside the neural network through which backpropagation occurs naturally. This is also in contrast to PINNs where the entire PDE, its boundary and its initial conditions are reproduced by the output of a neural network by adding them to the loss function. Importantly, the encoder learns an approximation of the inverse transform in a \emph{self-supervised} fashion that can be used to evaluate the hidden fields underlying unseen data without any further optimization. Moreover, the proposed framework is versatile as it allows for straightforward consideration of other domain-specific knowledge such as symmetries or constraints on the hidden field. In this work, we develop this idea for stationary and time-dependent PDEs on structured and unstructured grids and on noisy data using mesh-based and mesh-less numerical discretization methods.

\subsection{Novelties and features of BiPDEs}
A full PDE solver is implemented as a \textit{custom layer inside the architecture of semantic autoencoders} to solve inverse-PDE problems in a self-supervised fashion. Technically this is different than other works that implement a propagator decoder by manipulating activation functions or kernels/biases of convolutional layers, or those that feed the output of a neural network to a separate numerical solver such as in BN17 which requires the burden of considering the adjoint problem in order to compute partial differentiations. The novelties and features of this framework are summarized below:
\begin{enumerate}
\item \textbf{General discretizations.} We do not limit numerical discretization of differential equations to only finite differences that are emulated by convolution operations, our approach is more general and permits employing more sophisticated numerical schemes such as meshless discretizations. It is a more general framework that admits any existing discretization method directly in a decoder stage. 

\item \textbf{Introducing solver layers.} All the information about the PDE system is \emph{only} localized in a solver layer; \textit{i.e.} we do not inform the optimizer or the loss function with the adjoint PDE problem, or engineer regularizers or impose extra constraints on the kernels of convolutions, or define exotic activation functions as reviewed above. In other words, PDE solvers are treated as custom layers similar to convolution operations that are implemented in convolutional layers. An important aspect is the ability to employ any of the usual loss functions used in deep learning, for example we arbitrarily used mean absolute error or mean squared error in our examples.

\item \textbf{Blending meshless methods with trainable parameters.} Another unique proposal made in this work is the use of Radial Basis Function (RBF) based PDE solver layers as a natural choice to blend with deep neural networks. Contrary to other works, the neural network is not only used as an estimator for the unknown field but also it is tasked to optimize the shape parameters and interpolation points of the RBF scheme. In fact, our meshless decoder is not free of trainable parameters similar to reviewed works, instead shape parameters and seed locations are trainable parameters that define the RBF discretization, this is analogous to convolutional layers with trainable weights/biases that are used in machine learning domain. In fact this presents an example of neural networks complementing numerical discretization schemes. Choosing optimal shape parameters or seed locations is an open question in the field of RBF-based PDE solvers and here we show neural networks can be used to optimally define these discretization parameters. 

\item \textbf{Explicit/implicit schemes.} Most of the existing frameworks only accept explicit numerical discretizations in time, however our design naturally admits implicit methods as well. Using implicit methods allows taking bigger timesteps for stiff problems such as the diffusion problem, hence not only providing faster inverse-PDE solvers, but also present more robust/stable inverse PDE solvers.

\item \textbf{Data efficient.} Our design lowers the computational cost as a result of reusing classical numerical algorithms for PDEs during the learning process, which focuses provided data to infer the actual unknowns in the problem, \textit{i.e.} reduces the load of learning a discretization scheme from scratch.

\item \textbf{Physics informed.} Domain-specific knowledge about the unknown fields, such as symmetries or specialized basis functions, can be directly employed within our design.

\item \textbf{Inverse transform.} After training, the encoder can be used independently as a real-time estimator for unknown fields, \textit{i.e.} without further optimization. In other words, the network can be pre-trained and then used to infer unknown fields in real-time applications.

\end{enumerate}

\section{Blended inverse-PDE network (BiPDE-Net)}
The basic idea is to embed a numerical solver into a deep learning architecture to recover unknown functions in inverse-PDE problems, and all the information about the governing PDE system is only encoded inside the DNN architecture as a solver layer. In this section we describe our proposed architectures for inverse problems in one and two spatial dimensions.

\subsection{Deep neural networks (DNN)}
The simplest neural network is a single layer of perceptron that mathematically performs a linear operation followed by a nonlinear composition applied  to its input space,
\begin{align}
\mathcal{N} =\sigma\big( \mathbf{W}\mathbf{x}+\mathbf{b}\big),
\end{align}
where $\sigma$ is called the \textit{activation function}. Deep neural networks are multiple layers stacked together within some architecture. The simplest example is a set of layers connected in series without any recurrent loops, known as feedforward neural networks (FNN). In a densely connected FNN, the action of the network is simply the successive compositions of previous layer outputs with the next layers, \textit{i.e.},
\begin{align}
\mathcal{N}_l =\sigma\big( \mathbf{W}_l\mathcal{N}_{l-1}(\mathbf{x})+\mathbf{b}_l\big),
\end{align}
where $l$ indicates the index of a layer. This compositional nature of NNs is the basis of their vast potential as universal function estimators of any arbitrary function on the input space $\mathbf{x}$, see e.g. \cite{tikhomirov1991representation, cybenko1989approximation, csaji2001approximation}. Another important feature of NNs is that they can effectively express certain high dimensional problems with only a few layers, for example Darbon \textit{et al.} \cite{darbon2020overcoming} have used NNs to overcome the curse of dimensionality for some Hamilton-Jacobi PDE problems (also see \cite{han2018solving,Sirignano:2018aa}).

Most machine learning models are reducible to composition of simpler layers which allows for more abstract operations at a higher level. Common layers include dense layers as described above, convolutional layers in convolutional neural networks (CNNs) \cite{lecun1998gradient,krizhevsky2012imagenet}, Long-short term memory networks (LSTM) \cite{hochreiter1997long}, Dropout layers \cite{srivastava2014dropout} and many more. In the present work, we pay particular attention to CNNs owing to their ability to extract complicated spatial features from high dimensional input datasets. Furthermore, we define custom PDE solver layers as new member of the family of pre-existing layers by directly implementing numerical discretization schemes inside the architecture of deep neural networks.

\subsection{Custom solver layers}
A \textit{layer} is a high level abstraction that plays a central role in existing deep learning frameworks such as \texttt{TensorFlow}\footnote{For example see TensorFlow manual page at \href{https://www.tensorflow.org/guide/keras/custom_layers_and_models}{$\rm https://www.tensorflow.org/guide/keras/custom\_layers\_and\_models$}  } \cite{abadi2016tensorflow}, \texttt{Keras} API \cite{chollet2015keras}, \texttt{PyTorch} \cite{paszke2017automatic}, \textit{etc.} Each Layer encapsulates a state, \textit{i.e.} trainable parameters such as weights/biases, and a transformation of inputs to outputs. States in a layer could also be non-trainable parameters in which case they will be excluded from backpropagation during training. 

We implement different explicit or implicit numerical discretization methods as custom layers that transform an unknown field, initial data and boundary conditions to outputs in the solution space. Solver layers encapsulate numerical discretization schemes with trainable (\textit{e.g.} shape parameters and seeds in meshless methods) or non-trainable (\textit{e.g.} the finite difference methods) state parameters. Interestingly, solver layers with trainable parameters are new computational objects analogous to pre-existing convolutional layers with trainable kernel parameters. 

An important aspect of layer objects is that they can be composed with other layers in any order. Particularly, this offers an interesting approach for solving inverse problems given by systems of partial differential equations with several unknown fields that can be modeled with neural layers. We will explore this avenue in future work. In the remainder of this manuscript we will only focus on different inverse-PDE examples given by a single PDE equation and one unknown field.

\subsection{Blended neural network architectures}
\begin{figure}
\centering
\includegraphics[width=\linewidth]{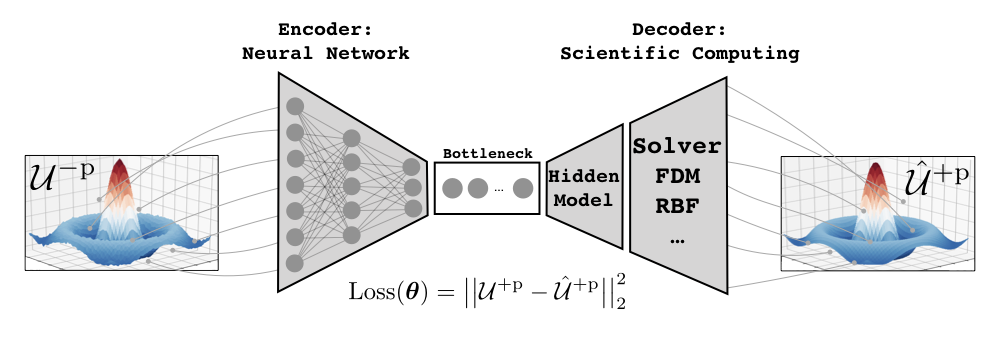}
\caption{Architecture of the BiPDE to infer unknown parameters of hidden fields. Here the loss function is the mean squared error between data and output of the autoencoder, however other choices for loss function may be used depending on the nature of data. }
\label{fig::autoBiPDE}
\end{figure}

BiPDE is a two-stage architecture, with the first stage responsible for learning the unknown coefficients and the second stage performing numerical operations as in traditional numerical solvers (see figure \ref{fig::autoBiPDE}). To achieve higher performance, it is essential to use GPU-parallelism. We leverage the capability provided by the publicly available library \texttt{TensorFlow} \cite{abadi2016tensorflow} by implementing our PDE-solver as a \textit{custom layer} into our network using the \texttt{Keras} API \cite{chollet2015keras}. Details of this includes vectorized operations to build the linear system associated by the PDE discretization.

We propose a semantic autoencoder architecture as proposed by Aragon-Calvo (2019) \cite{aragon2019self} with hidden parameters being represented at the bottleneck of the autoencoder. Figure \ref{fig::autoBiPDE} illustrates the architecture for the proposed semantic autoencoder. Depending on static or time dependent nature of the governing PDE, one may train this network over pairs of input-output solutions that are shifted $\rm p$ steps in time, such that for a static PDE we have $\rm p=0$ while dynamic PDEs correspond to $\rm p\ge 1$. We call this parameter the \textit{shift parameter}, which will control the accuracy of the method (\emph{cf.} see section \ref{sec::meshfree}).

An important aspect is that the input to BiPDE is the solution data itself. In other words the neural network in a BiPDE is learning the \textit{inverse transform},
\begin{align}
\rm \mathcal{NN}:~ \{u\}\rightarrow \textrm{hidden\ field},
\end{align}
where $\{u\}$ indicates an ensemble of solutions, \textit{e.g.} solutions obtained with different boundary conditions or with different hidden fields. Note that in other competing methods such as PINNs the input is sanctioned to be the coordinates in order for automatic differentiation to compute spatial and temporal derivatives; as a consequence PINNs can only be viewed as \textit{surrogates} for the solution of the differential problem defined on the space of coordinates. However, we emphasize that semantic autoencoders are capable to approximate the inverse transformation from the space of solutions to the space of hidden fields, a feature that we exploit in section \ref{sec::inverse}.

Essentially different numerical schemes can be implemented in the decoder stage. We will blend examples of both mesh-based and mesh-less numerical discretizations and present numerical results and comparisons with PINNs. We will show how BiPDEs can handle data on unstructured grids and data with added noise. In section \ref{sec::meshbased}, we demonstrate performance of mesh-based BiPDEs on inverse problems in two spatial dimensions by using a finite difference discretization and Zernike expansion of the non-homogeneous hidden field, we will consider both stationary and dynamic PDE problems in this section. Then in section \ref{sec::meshfree}, we develop a mesh-less BiPDE and consider a dynamic nonlinear inverse partial differential problem. 

\section{Mesh-based BiPDE: Finite Differences}\label{sec::meshbased}
We consider a variable coefficient Poisson problem in one and two spatial dimensions as well as the one dimensional nonlinear Burger's equation as an example of a nonlinear dynamic PDE problem with a scalar unknown parameter. 

\subsection{Stationary Poisson problem}
We consider the governing equation for diffusion dominated processes in heterogeneous media:
\begin{align}
&\nabla\cdot \big(D(\mathbf{x}) \nabla u\big)=-f(\mathbf{x}), &\mathbf{x}\in \Omega \label{eq::Poisson}\\
&u(\mathbf{x})=u_0(\mathbf{x}), &\mathbf{x}\in\partial \Omega
\end{align}
Here we consider a rectangular domain with Dirichlet boundary conditions.

\textbf{Discretization.} In our architecture, we use the standard 5-point stencil finite difference discretization of the Poisson equation in the solver layer, \textit{i.e.}
\begin{align*}
&\frac{D_{i-1/2, j} u_{i-1, j} - (D_{i-1/2,j} + D_{i+1/2,j})u_{i,j} + D_{i+1/2,j}u_{i+1,j}}{\Delta x^2}+\\
&\frac{D_{i, j-1/2} u_{i, j-1} - (D_{i,j-1/2} + D_{i,j+1/2})u_{i,j} + D_{i,j+1/2}u_{i,j+1}}{\Delta y^2} + f_{i,j}=0,
\end{align*}
and we use the linear algebra solver implemented in \texttt{TensorFlow} to solve for the solution field, \textit{i.e.} we used \texttt{tf.linalg.solve} method that is a dense linear system solver. Of course, this can be improved by implementing a sparse linear solver. 

\textbf{Hidden Model.} We decompose the hidden field into a finite number of eigenfunctions and search for their optimal coefficients. This is also advantageous from a physics point of view, because domain's knowledge of hidden fields can be naturally formulated in terms of basis functions into this framework. One such family of series expansions are the moment-based methods that have been largely exploited in image reconstruction \cite{khotanzad1990invariant, belkasim1991pattern, prokop1992survey, bailey1996orthogonal}. In particular, Zernike moments \cite{von1934beugungstheorie} provide a linearly independent set of polynomials defined on the unit circle/sphere in two/three spatial dimensions. Zernike moments are well-suited for such a task and are commonly used for representing optical aberration in astronomy and atmospheric sciences \cite{ragazzoni2000adaptive}, for image reconstruction and for enhanced ultrasound focusing in biomedical imaging \cite{dong2008zernike,markelj2012review,kaye2012application}.   

Zernike moments are advantageous over regular moments in that they intrinsically provide rotational invariance, higher accuracy for irregular patterns, and are orthogonal, which reduces information redundancy in the different coefficients. Zernike polynomials capture deviations from zero mean as a function of radius and azimuthal angle. Furthermore, the complete set of orthogonal bases provided by Zernike moments can be obtained with lower computational precision from input data, which enhances the robustness of the reconstruction procedure. 

Odd and even Zernike polynomials are given as a function of the azimuthal angle $\theta$ and the radial distance $\rho$ between $0$ and $1$ measured from the center of image,
\begin{align*}
&\begin{bmatrix} Z_{nm}^o(\rho, \theta) \\ Z_{nm}^e(\rho, \theta) \end{bmatrix}=R_{nm}(\rho) \begin{bmatrix} \sin(m\theta)\\ \cos(m\theta)\end{bmatrix},
\end{align*}
with
\begin{align*}
R_{nm}(\rho)&=\begin{cases}
\sum_{l=0}^{(n-\vert m\vert )/2}\frac{(-1)^l (n-l)!}{l![(n+\vert m\vert )/2 -l]! [(n-\vert m\vert )/2 -l]!}\rho^{n-2l}  & \textrm{for $n-m$ even,} \\
 0  &  \textrm{for $n-m$ odd,}
\end{cases}
\end{align*}
where $n$ and $m$ are integers with $n\ge \vert m\vert$. A list of radial components is given in table \ref{tab::I} (from \cite{weisstein2002zernike}). For an extensive list of Zernike polynomials in both two and three spatial dimensions, we refer the interested reader to \cite{mathar2008zernike}.
\begin{table}
\centering
\resizebox{1 \textwidth}{!}{
\begin{tabular}{SSSSSS} \toprule 
    {$\textbf{n}$} & {$ \vert \textbf{m}\vert$} & {$\textbf{R}_{nm}$}                 & {$\textbf{Z}_{nm}^o$}                                  & {$\textbf{Z}_{nm}^e$}                                 & {$\rm \textbf{Aberration/Pattern}$}\\ \midrule \midrule
    0  & 0         & 1                                       & {$0$}                                                & 1                                                    & {$\rm Piston$} \\  \midrule
    1  & 1         & {$\rho$}                             & {$\rho\sin(\theta)$}                          &  {$\rho\cos(\theta)$}                       & {$\rm Tilt $}\\  \midrule
    2  & 0         &  {$2\rho^2 - 1$}                 &  {$0$}                                              &  {$2\rho^2 - 1$}                               & {$\rm Defocus $}\\
        & 2         & {$\rho^2$}                         &  {$\rho^2\sin(2\theta)$}                   &  {$\rho^2\cos(2\theta)$}                   & {$\rm Oblique/Vertical\ Astigmatism $} \\ \midrule
    3  & 1         & {$3\rho^3 - 2\rho$}            &  {$(3\rho^3 - 2\rho)\sin(\theta)$}     &  {$(3\rho^3 - 2\rho)\cos(\theta)$}      & {$\rm Vertical/Horizontal\ Coma $}\\
        & 3         & {$\rho^3$}                         & {$\rho^3\sin(3\theta)$}                     & {$\rho^3\cos(3\theta)$}                     & {$\rm Vertical/Oblique\ Trefoil $}\\  \midrule
    4  & 0         & {$6\rho^4 - 6\rho^2 + 1 $} &  {$0$}                                               & {$6\rho^4 - 6\rho^2 + 1 $}                  &  {$\rm Primary\ Spherical $} \\
        & 2         & {$4\rho^4 - 3\rho^2$}        &  {$(4\rho^4 - 3\rho^2)\sin(2\theta)$} &  {$(4\rho^4 - 3\rho^2)\cos(2\theta)$}    & {$\rm Oblique/Vertical\ Secondary\ Astigmatism $} \\ 
        & 4         & {$\rho^4$}                         &  {$\rho^4\sin(4\theta)$}                     &  {$\rho^4\cos(4\theta)$}                        & {$\rm Oblique/Vertical\ Quadrafoil $}\\ \bottomrule 
\end{tabular}
}
     \caption{First $15$ odd and even Zernike polynomials according to Noll's nomenclature. Here, the ordering is determined by ordering polynomial with lower radial order first, cf. \cite{wyant1992basic}.\label{tab::I}}
\end{table}

Furthermore, each Zernike moment is defined by projection of the hidden field $f(x,y)$ on the orthogonal basis,
\begin{align*}
& \begin{bmatrix} A_{nm}\\ B_{nm}\end{bmatrix}=\frac{n+1}{\epsilon_{mn}^2\pi}\int_x \int_y f(x,y) \begin{bmatrix} Z_{nm}^o(x,y) \\  Z_{nm}^e(x,y) \end{bmatrix} dx dy, \quad x^2 + y^2\le 1,
\end{align*}
where for $m=0,~ n\neq 0$ we defined $\epsilon_{0n}=1/\sqrt{2}$ and $\epsilon_{mn}=1$ otherwise. Finally, superposition of these moments expands the hidden field in terms of Zernike moments:
\begin{align}
\hat{f}(x,y)=\sum_{n=0}^{N_{max}}\sum_{\vert m\vert=0}^{n}\big[  A_{nm}Z_{nm}^o (r,\theta) + B_{nm}Z_{nm}^e (r,\theta) \big]. \label{eq::ZExp}
\end{align}

\begin{figure}
\centering
\includegraphics[width=\linewidth]{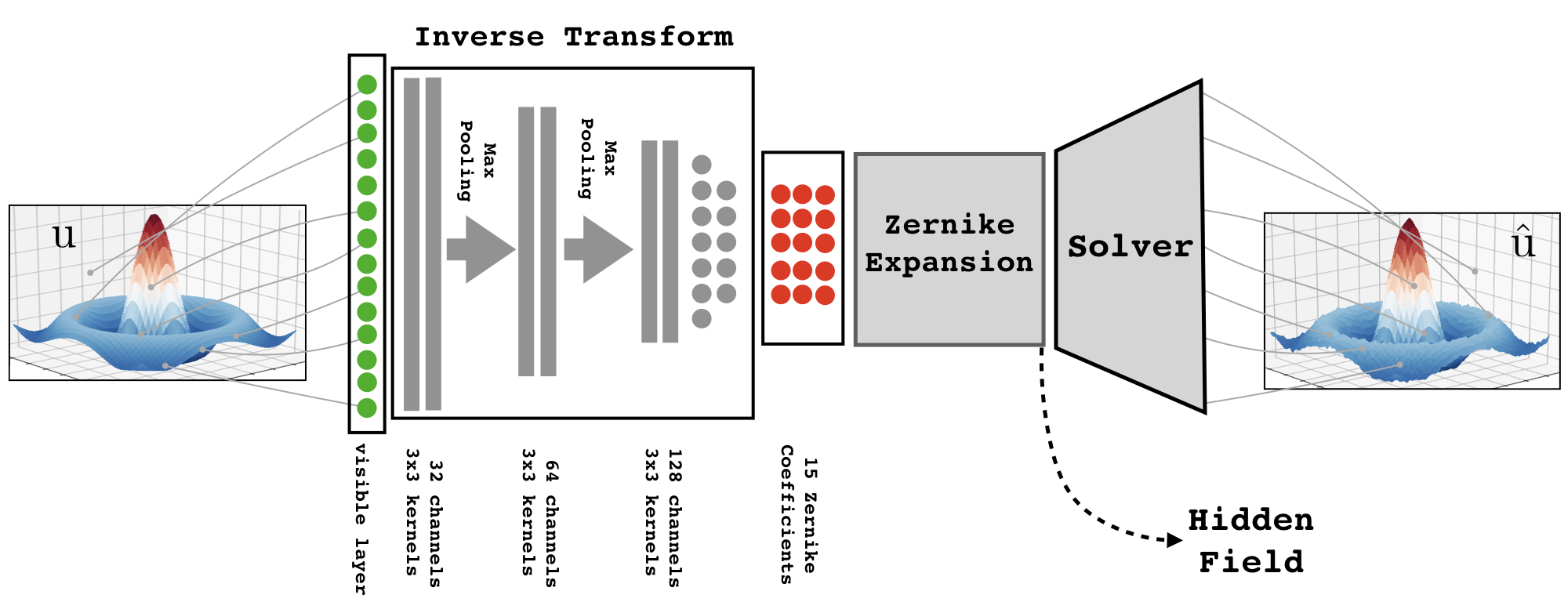}
\caption{Architecture of the semantic autoencoder to infer hidden fields. Zernike moments are discovered at the bottleneck of the architecture.}
\label{fig::autoarch}
\end{figure}

In order to identify the coefficients in the Zernike expansion \eqref{eq::ZExp} for hidden fields, we use a semantic autoencoder architecture with Zernike moments being represented by the code at the bottleneck of the autoencoder. Figure \ref{fig::autoarch} illustrates the architecture for the proposed semantic autoencoder.

\textbf{Architecture.} Even though a shallow neural network with as few neurons as the number of considered Zernike terms suffices to estimate values of the unknown Zernike moments in each of the problems considered in this section, however we will use a deep convolutional neural network (detailed below) in order to achieve our ultimate goal of approximating the inverse transform for the Poisson problem in a broad range of diffusion coefficient fields. Therefore we design one deep neural network and uniformly apply it to several problems in this section.

In the training process of a CNN, the kernels are trained at each layer such that several feature maps are extracted at each layer from input data. The CNN is composed of 3 convolutional blocks with $32,~ 64,~ 128$ channels respectively and kernel size $3\times 3$. Moreover, we use the \texttt{MaxPooling} filter with kernel size $(2,2)$ after each convolutional block to downsample the feature maps by calculating the maximum values of each patch within these maps. We use the \texttt{ReLU} activation function \cite{hahnloser2000digital}, \textit{i.e.} a piecewise linear function that only outputs positive values: ${\rm ReLU}(x)=\max(0,x)$, in the convolutional layers followed by a \texttt{Sigmoid} activation in dense layers and a scaled \texttt{Sigmoid} activation at the final layer,
\begin{align}
\tilde{\sigma}(x)&=D_{\min} + (D_{\max} - D_{\min})\sigma (x)\label{eq::activ},
\end{align}
such that the actual values of the diffusion coefficient are within the range $(D_{\min}, D_{\max})$, known from domain specific knowledge. After each dense layer, we apply \texttt{Dropout} layers with a rate of $0.2$ to prevent overfitting \cite{hinton2012improving,srivastava2014dropout} (a feature that is most useful in estimating the inverse transform operator) and avoid low quality local minima during training.
 
\subsubsection{Test cases.} \label{sec::tests}
\textbf{Case I. A tilted plane.}
In the first example we consider a linear diffusion model given by
\begin{align*}
&D(x,y)=\sqrt{2} + 0.1(y-x)
\end{align*}
where the boundary condition function $u_{BC}$ and the source field $f$ are given by
\begin{align*}
&u_{BC}(x,y)=0.01\cos(\pi x)\cos(\pi y)  \qquad \textrm{and} \qquad f(x,y)=\sin(\pi x)\cos(\pi y)
\end{align*}
\begin{figure}
\centering
\subfigure[Comparison of learned (left) versus true diffusion coefficient (right).]{\includegraphics[width=\linewidth]{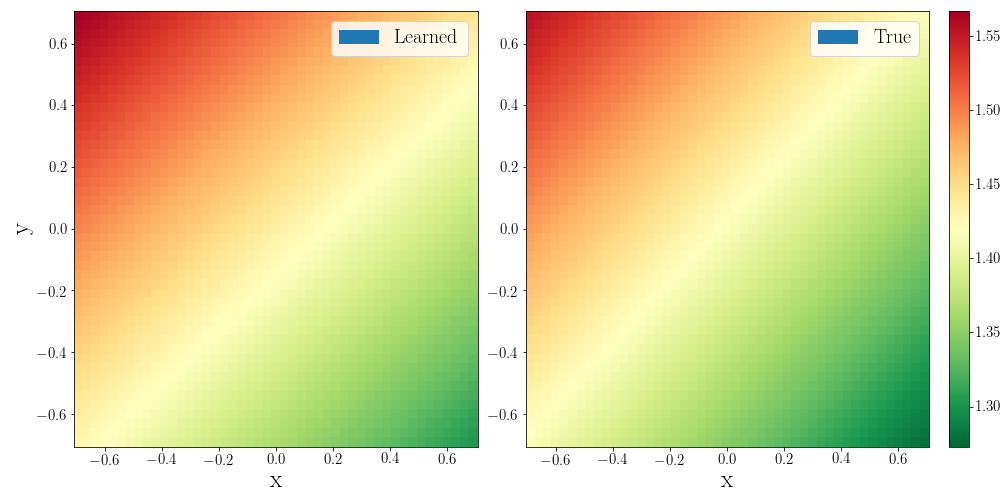} \label{subfig::b}}\quad\quad
\subfigure[Learned solution.]{ \includegraphics[width=0.45\linewidth]{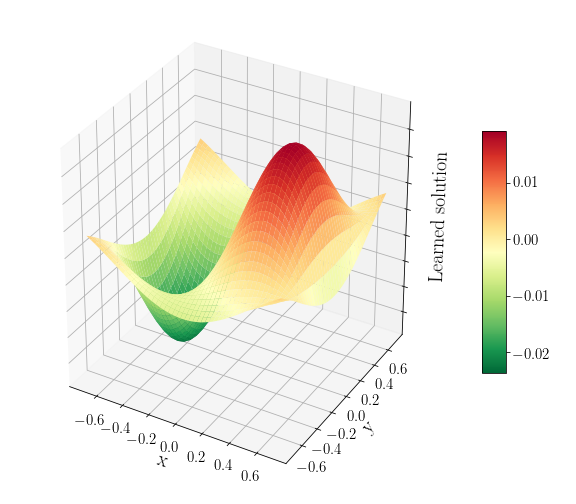} }\quad\quad
\subfigure[True solution.]{ \includegraphics[width=0.45\linewidth]{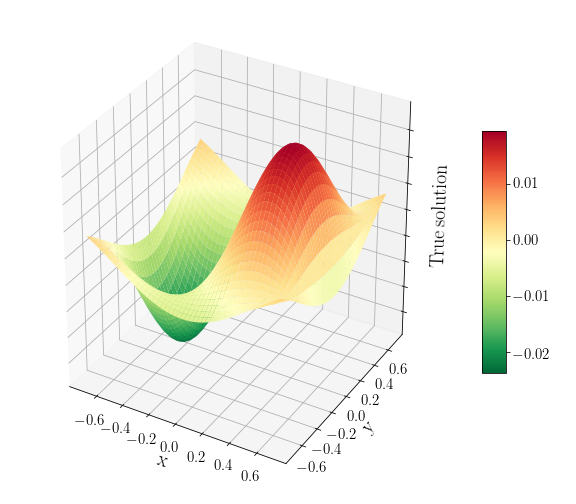} }\quad\quad
\subfigure[Error in learned solution $u-\hat{u}$.]{\includegraphics[width=0.45\linewidth]{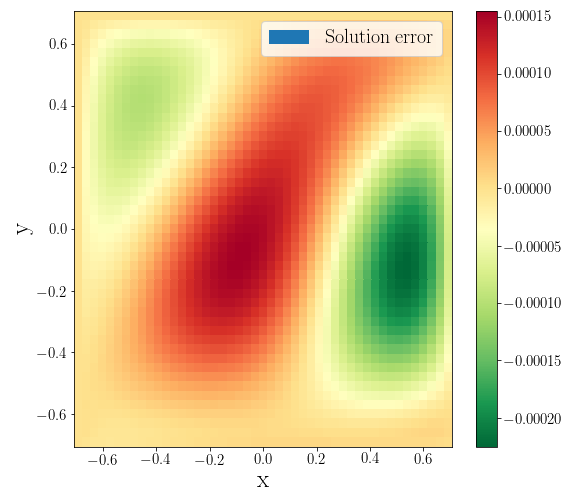} }\quad\quad
\subfigure[Error in learned diffusion coefficient.]{\includegraphics[width=0.45\linewidth]{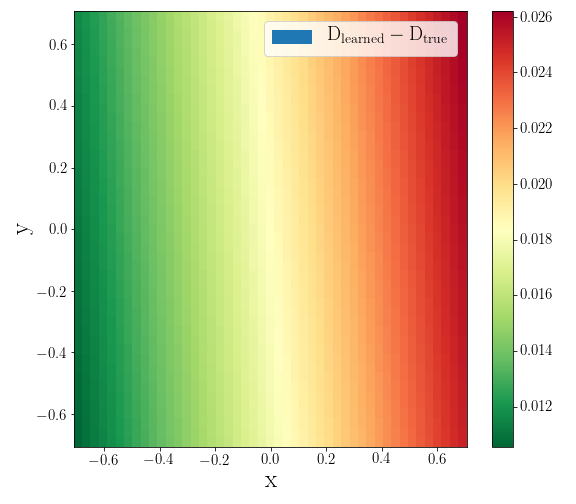} }\quad\quad
\caption{Results for the two dimensional tilted plane (case I).}
\label{fig:example2}
\end{figure}
In this experiment we only use a \emph{single solution field} for training. Even though in our experiments the method succeeded to approximate the hidden field even with a \emph{single grid point} to compute the loss function, here we consider all the grid points in the domain to obtain improved accuracy in the results. We trained the network for $\rm 30$ epochs using an \texttt{Adam} optimizer \cite{kingma2014adam} that takes $170$ seconds on a Tesla T4 GPU available on a free Google Colaboratory account\footnote{\href{https://colab.research.google.com/}{https://colab.research.google.com/}}.  Figure \ref{fig:example2} depicts the results obtained by the proposed scheme. The diffusion map is discovered with a maximum relative error of only $2\%$, while the error in the solution field is $1\%$. It is noteworthy to mention that the accuracy of the results in this architecture are influenced by the accuracy of the discretizations used in the solver layer. While we used a second-order accurate finite difference discretization, it is possible to improve these results by using higher order discretizations instead. We leave such optimizations to future work.
\begin{figure}
\centering
\includegraphics[width=0.55\linewidth]{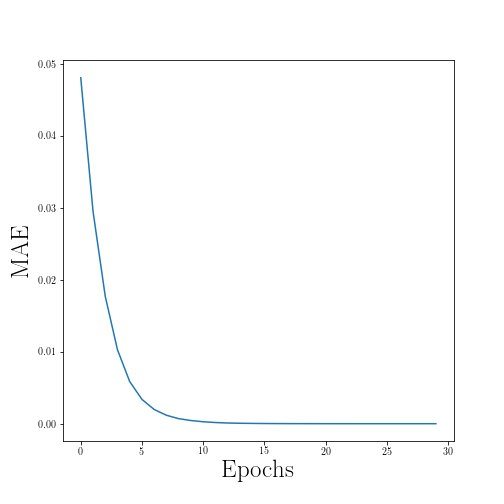} 
\caption{Mean absolute error in solution vs. epochs for the two dimensional tilted plane (case I).}
\label{fig:example2loss}
\end{figure}

\begin{table}[ht]
\centering
\resizebox{1 \textwidth}{!}{
\begin{tabular}{| c | c | c | c || c |  c || c | c ||| c| c | c |c|}\toprule \hline
T &  \# params &  \texttt{C(32)} & \texttt{C(32)} & \texttt{C(64)} & \texttt{C(64)} & \texttt{C(128)} & \texttt{C(128)} & \texttt{D(64)} & \texttt{D(32)} & $\rm MAE_D$ & $\rm L^\infty_D$ \\ \hline\hline
1 & $1,468,323$ & Y & Y & Y & Y & Y & Y & Y & Y & $0.0144207$ & $0.0294252$\\
2 & $1,459,075$ & Y & - & Y & Y & Y & Y & Y & Y & $0.0193128$ & $0.0267854$\\
3 & $1,422,147$ & Y & - & Y & - & Y & Y & Y & Y & $0.0226252$ & $0.0527432$\\
4 & $1,274,563$ & Y & - & Y & - & Y & - & Y & Y & $0.0199361$ & $0.0272122$ \\
5 & $682,627$   & Y & - & Y & - & Y & - & - & Y & $0.0141946$ & $0.0243868$ \\
6 & $313,859$   & Y & - & Y & - & - & - & - & Y & $0.0301841$ & $0.0544990$ \\
7 & $46,467$    & Y & - & Y & - & - & - & - & - & $0.0190432$ & $0.0264254$ \\
8 & $6,915$     & - & - & - & - & - & - & - & - & $0.0183808$ & $0.0267156$ \\
\hline\bottomrule
\end{tabular}
}
\caption{Influence of architecture of the decoder stage on mean absolute error $\rm MAE_D \equiv \rm \sum \vert D(\mathbf{x}) - \hat{D}(\mathbf{x})\vert/N$ and maximum error $\rm L^\infty_D$ in the discovered hidden field in case I. Double vertical lines correspond to \texttt{MaxPooling2D()} layers and triple vertical lines correspond to \texttt{Flatten()} layer. \texttt{C(o)} and \texttt{D(o)} stand for \texttt{conv2D(filters)} and \texttt{Dense(neurons)} layers respectively. There are $3$ neurons at the bottleneck not shown in the table. }
\label{tab::arch_field}
\end{table}

\textbf{Influence of architecture.} Table \ref{tab::arch_field} tabulates the mean absolute error in the discovered tilted plane diffusion coefficient for different architectures of the encoder stage. No significant improvement is observed for deeper or shallower encoder network for the example considered here.

\begin{figure}
\centering
\subfigure[Learned diffusion.]{ \includegraphics[width=0.45\linewidth]{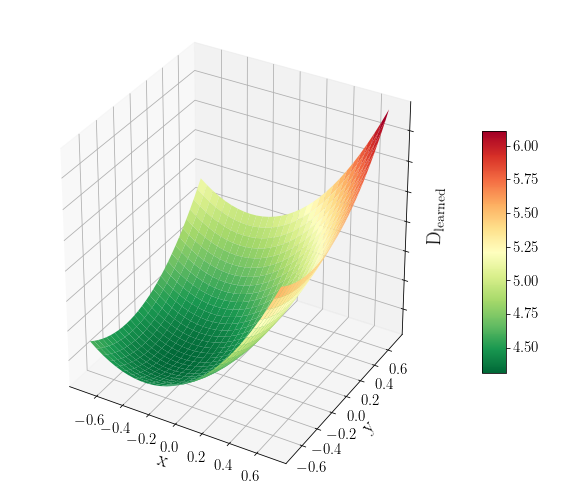} }\quad\quad
\subfigure[True diffusion.]{ \includegraphics[width=0.45\linewidth]{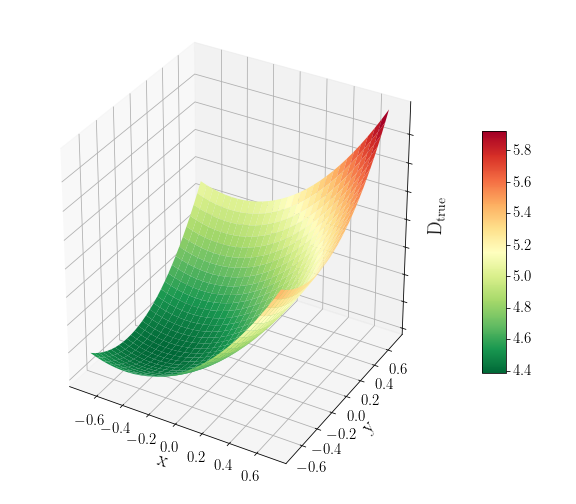} }\quad\quad
\subfigure[Learned solution.]{ \includegraphics[width=0.45\linewidth]{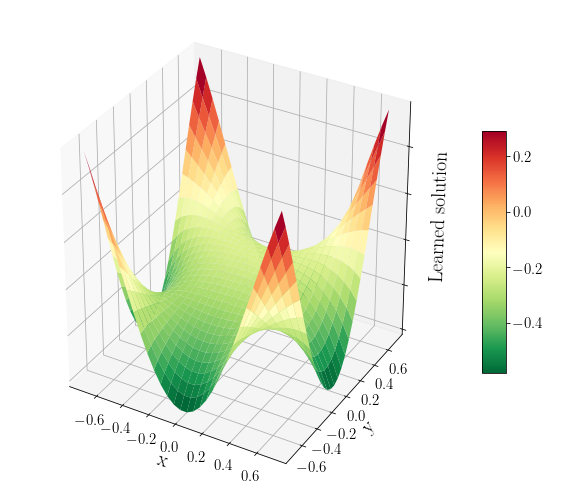} }\quad\quad
\subfigure[True solution.]{ \includegraphics[width=0.45\linewidth]{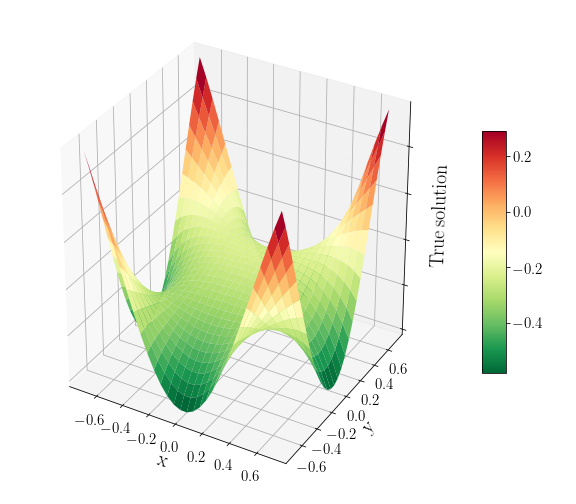} }\quad\quad
\subfigure[Error in learned solution $u-\hat{u}$.]{\includegraphics[width=0.45\linewidth]{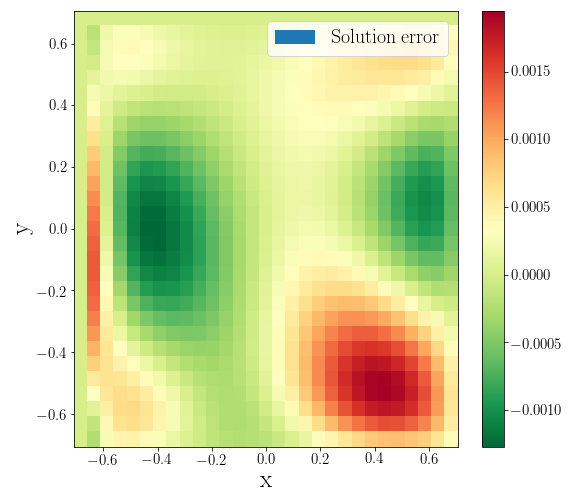} }\quad\quad
\subfigure[Error in learned diffusion coefficient.]{\includegraphics[width=0.45\linewidth]{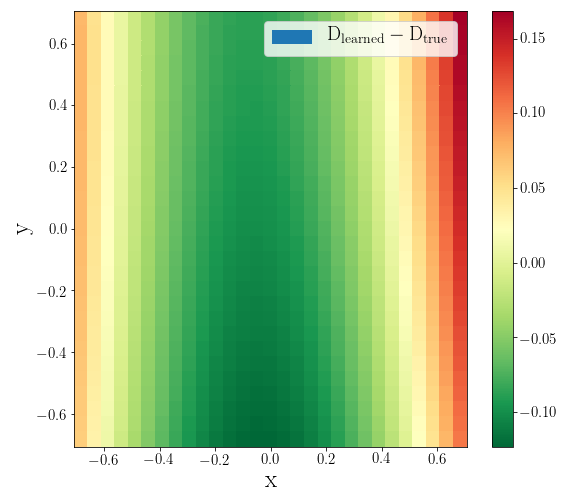} }\quad\quad
\caption{Results in the two dimensional parabolic case.}
\label{fig:example3}
\end{figure}

\textbf{Case II. superimposed Zernike polynomials.}
We consider a more complicated hidden diffusion field given by 
\begin{align*}
&D(x,y)= 4 + a_0 + 2a_1x + 2 a_2 y + \sqrt{3}a_3 (2 x^2 + 2y^2 - 1).
\end{align*}
The boundary condition function $u_{BC}$ and the source field $f$ are given by
\begin{align*}
&u_{BC}(x,y)=\cos(\pi x)\cos(\pi y)\qquad \textrm{and} \qquad f(x,y)=x+ y.
\end{align*}
Figure \ref{fig:example3} illustrates the performance of the proposed Zernike-based network using a mean absolute error measure for the loss function. We trained the network for $\rm 100$ epochs using an \texttt{Adam} optimizer \cite{kingma2014adam}.

\begin{figure}
\centering
\subfigure[Learned diffusion.]{ \includegraphics[width=0.45\linewidth]{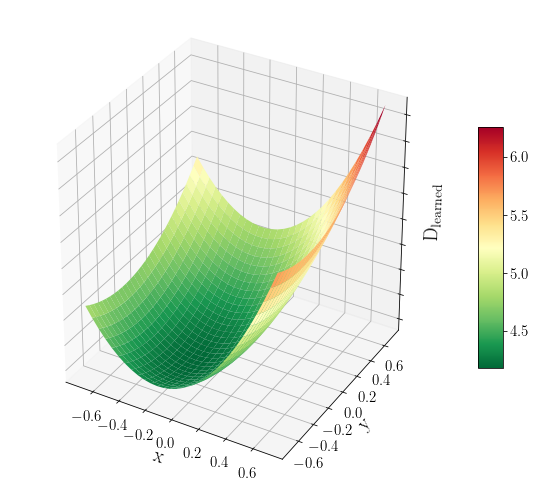} }\quad\quad
\subfigure[True diffusion.]{ \includegraphics[width=0.45\linewidth]{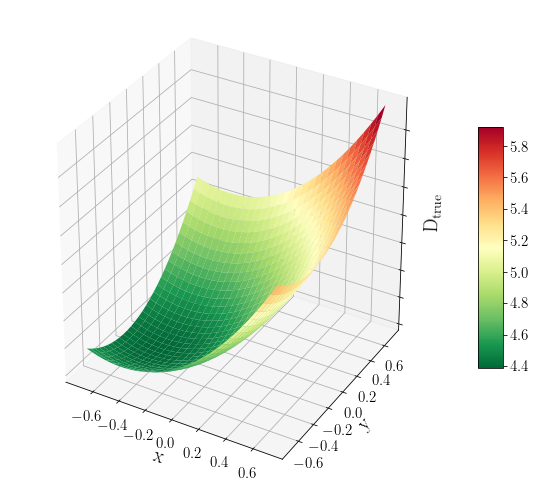} }\quad\quad
\subfigure[Learned solution.]{ \includegraphics[width=0.45\linewidth]{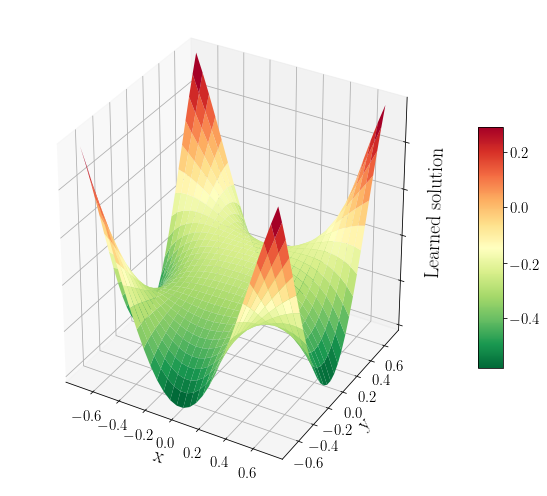} }\quad\quad
\subfigure[Noisy input solution.]{ \includegraphics[width=0.45\linewidth]{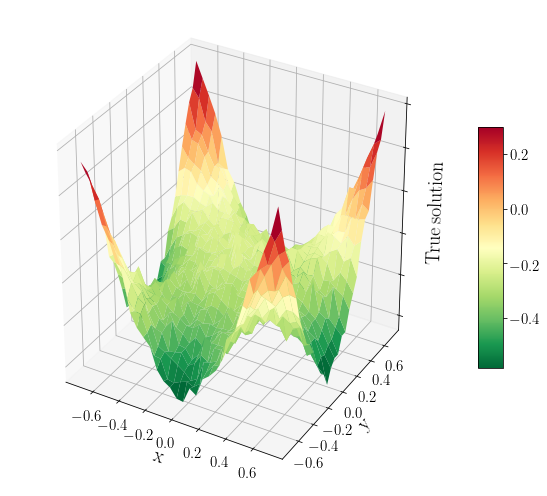} }\quad\quad
\subfigure[Error in learned solution $u-\hat{u}$.]{\includegraphics[width=0.45\linewidth]{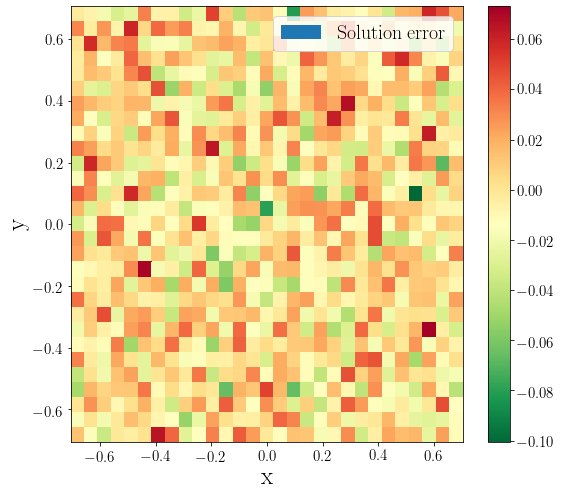} }\quad\quad
\subfigure[Error in learned diffusion coefficient.]{\includegraphics[width=0.45\linewidth]{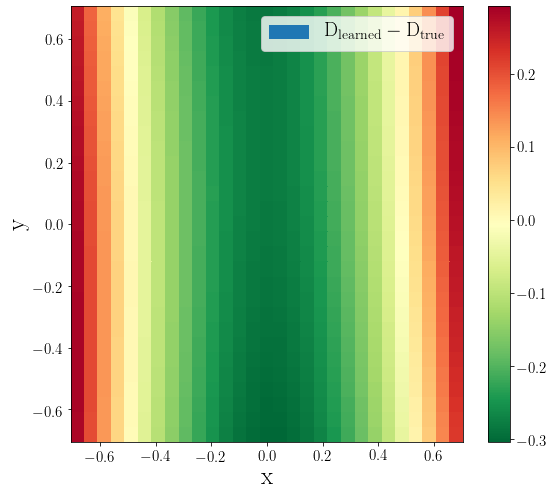} }\quad\quad
\caption{Results in the two dimensional case with added noise. After 300 epochs the network discovers the hidden diffusion field with a maximum relative error of $5\%$. Interestingly the learned solution is resilient to added noise and the network approximates a noise-free solution.}
\label{fig:performance}
\end{figure}

\textbf{Resilience to noise.} We also assess the performance of our scheme on noisy datasets. We consider a zero-mean Gaussian noise with standard deviation $0.025$ superimposed on the solution field. Figure \ref{fig:performance} depicts the solution learned from a noisy input image. The network succeeds in discovering the diffusion field with comparable accuracy as in the noise-free case. Note that this architecture naturally removes the added noise from the learned solution, a feature that is similar to applying a low-pass filter on noisy images.

\begin{figure}
\centering
\subfigure[$\rm L_1$ loss vs. epoch for case II without added noise.]{\includegraphics[width=0.45\linewidth]{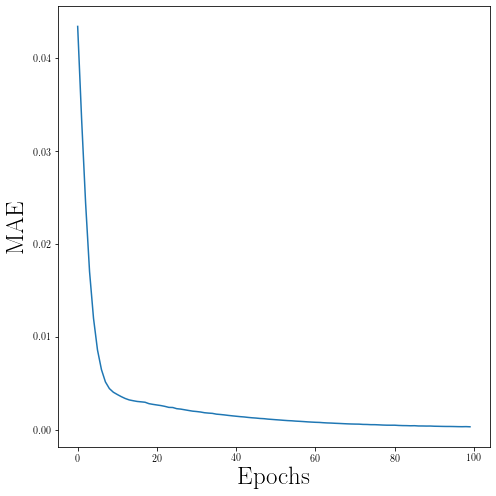} \label{subfig::losspara}}\quad\quad
\subfigure[$\rm L_2$ loss vs. epoch for case II with added noise.]{\includegraphics[width=0.45\linewidth]{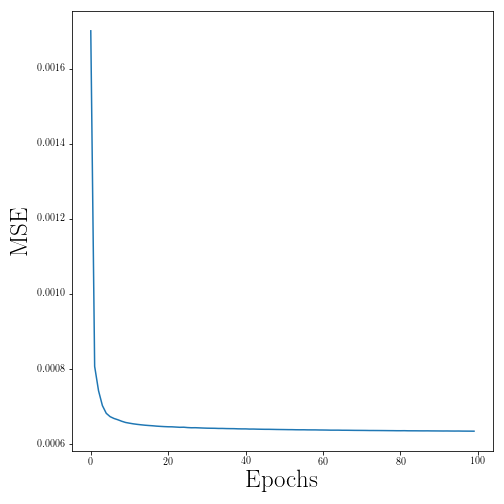} \label{subfig::lossparaNose}}\quad\quad
\subfigure[$\rm L_2$ loss vs. epoch for 1D inverse transform.]{\includegraphics[width=0.45\linewidth]{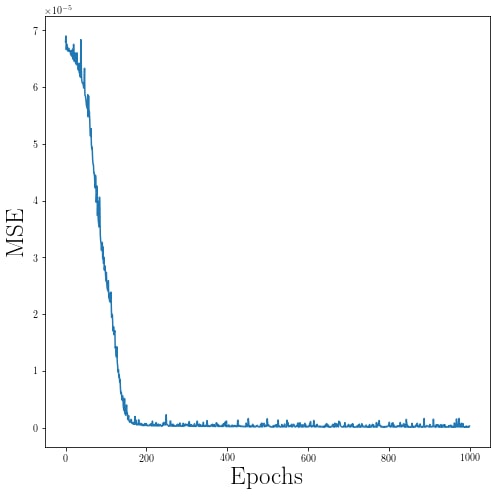} \label{subfig::lossinv1D}}\quad\quad
\subfigure[$\rm L_2$ loss vs. epoch for 2D inverse transform.]{\includegraphics[width=0.45\linewidth]{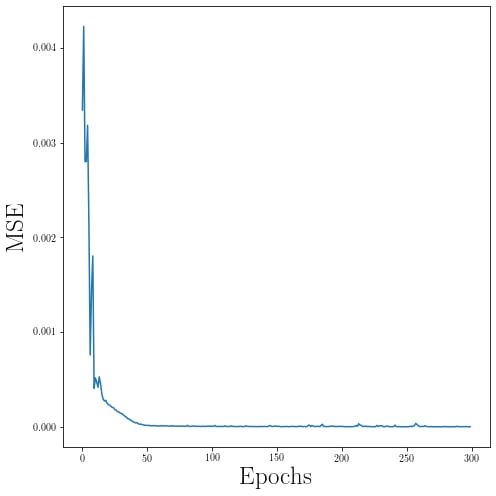} \label{subfig::lossinv2D}}\quad\quad
\caption{Mean absolute/square error vs. epochs for (top panel) the two dimensional parabolic experiment (case II) with and without added Gaussian noise of section \ref{sec::tests}, and (bottom panel) the inverse transform for 1D and 2D experiments of section \ref{sec::inverse}.}
\label{fig:example3loss}
\end{figure}

\subsubsection{Learning the inverse transform}\label{sec::inverse}
In the previous sections, we have applied BiPDE to find the variable diffusion coefficient from a single input image. Another interesting feature of the proposed semantic autoencoder architecture is its ability to train neural networks in order to discover the inverse transform for the underlying hidden fields \textit{in a self-supervised fashion}. In this scenario, the trained encoder learns the inverse transform function that approximates the hidden parameters given a solution field to its input. Note that even though the same task could be accomplished by supervised learning of the hidden fields, \emph{i.e.} by explicitly defining loss on the hidden fields without considering the governing equations, BiPDEs substitute the data labels with a governing PDE and offer comparable prediction accuracy. In this section we train BiPDEs over ensembles of solution fields to estimate hidden Zernike moments of diffusion coefficients underlying unseen data.

\textbf{One dimensional inverse transform}

\begin{figure}
\centering
\subfigure[Regression quality is $\rm R^2=0.9906891$.]{ \includegraphics[width=0.45\linewidth]{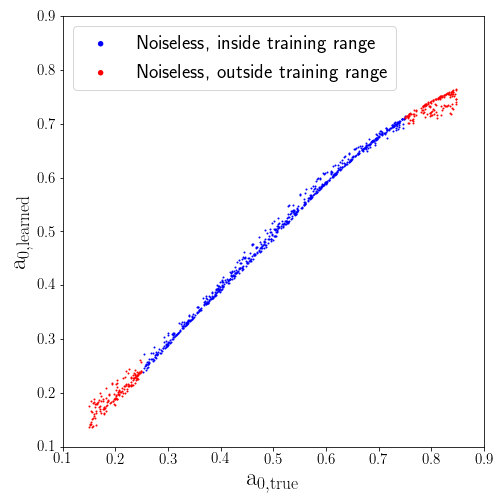} } 
\subfigure[Regression quality is $\rm R^2=0.9953392$.]{ \includegraphics[width=0.45\linewidth]{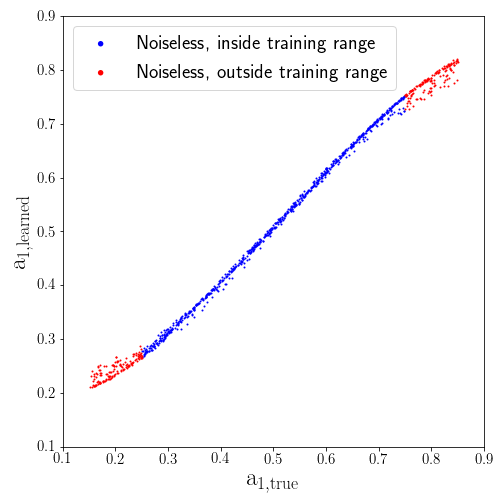} }
\subfigure[Regression quality is $\rm R^2=0.9796781$.]{ \includegraphics[width=0.45\linewidth]{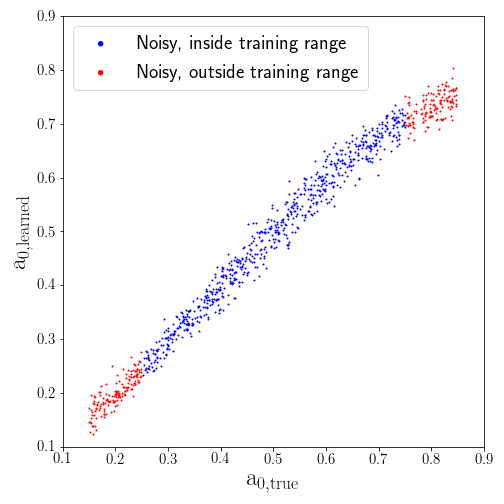} } 
\subfigure[Regression quality is $\rm R^2=0.9834912$.]{ \includegraphics[width=0.45\linewidth]{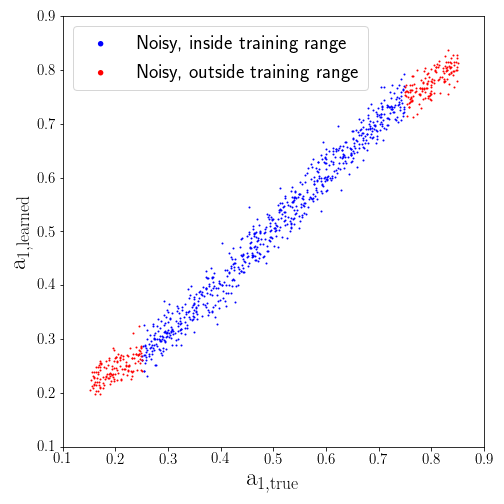} }
\caption{(Top, bottom) panel shows performance of BiPDE over $1000$ randomly chosen one-dimensional images with $\rm N_x=160$ grid points after $1000$ epochs (with,without) added zero-mean Gaussian noise with standard deviation $0.025$ to the test sample. The hidden diffusion coefficient is $D(x)=1 + a_0 + a_1 x$. In each case the $\rm R^2$ coefficient is reported for the blue data points, where unknown parameters fall within the training range $[0.25,0.75]$. Red data points show predictions outside of training range. Network has $20,222$ trainable parameters, and training takes $\sim 2$ seconds per epoch on a Tesla T4 GPU available on a free Google Colaboratory account.}
\label{fig:1d_inverse}
\end{figure}
We build a one dimensional semantic autoencoder using 3 layers with $100,~ 40$, and $2$ neurons respectively. We used the \texttt{ReLU} activation function for the first two layers and a \texttt{Sigmoid} activation function for the last layer representing the hidden parameters. A linear solver is then stacked with this encoder that uses the second order accurate finite difference discretization, \textit{i.e.}
\begin{align*}
&\frac{D_{i-1/2} u_{i-1} - (D_{i-1/2} + D_{i+1/2})u_{i} + D_{i+1/2}u_{i+1} }{\Delta x^2} + f_{i}=0, & D_{i+1/2}=\frac{D_i + D_{i+1}}{2}
\end{align*}
However, the diffusion map is internally reconstructed using the hidden parameters before feeding the output of the encoder to the solver. As a test problem, we consider the one dimensional Poisson problem with a generic linear form for the diffusion coefficient,
\begin{align*}
D(x)=1 + a_0 + a_1 x.
\end{align*}
We consider identical left and right Dirichlet boundary conditions of 0.2 for all images and let the source term be $f(x)=\sin(\pi x)$. We consider random diffusion coefficients $a_0$ and $a_1$ with a uniform distribution in $[0.25, 0.75]$ and we generate $1000$ solutions over the domain $x\in [-1,1]$. We train BiPDE over $900$ images from this dataset and validate its performance over the remaining $100$ images using a mean squared error loss function for $1000$ epochs. Each image is generated on a uniform grid with $\rm N_x=160$ grid points. We used a batch size of $100$ in these experiments using the \texttt{Adam} optimizer. Figure \ref{subfig::lossinv1D} shows loss versus epochs in this experiment. Figure \ref{fig:1d_inverse} compares learned and true coefficients over two independent test samples containing $1000$ solutions, with and without a zero-mean Gaussian noise with standard deviation $0.025$, \emph{i.e.} amounting to $\sim 13\%$ added noise over the images.  

In figure \ref{fig:1d_inverse}, we expanded the range of unknown parameters $a_0,a_1\in [0.15, 0.85]$ in our test sample to assess performance of trained encoder over unseen data that are outside the range of training data (as a measure of generalizability of BiPDEs). In this figure blue points correspond to new images whose true unknown parameters fall inside the training range, and red data points correspond to those outside the training range. We observe that the encoder is able to predict the unknown parameters even outside of its training range, although its accuracy gradually diminishes far away from the training range. Note that the predicted values for $a_0$ and $a_1$ exhibit a systematic error towards the lower and upper bounds of the \texttt{Sigmoid} activation function, indicative of the influence of the \texttt{Sigmoid} activation function used in the last layer. This emphasizes the significance of properly designing activation functions at the bottleneck.

Using the $\rm R^2$ statistical coefficient as a measure of accuracy for the trained encoder, we assess effects of sample size and grid points on the performance of BiPDEs and report the results in table \ref{tab::tab::R2coeff1D}.
\begin{enumerate}
\item \textit{Effect of sample size:} First, we fix number of grid points and vary sample size. We find that increasing sample size improves accuracy of the predictions in the case of clean data, however in the case of noisy data the accuracy does not show significant improvement by enlarging sample size.

\item \textit{Effect of grid points:} Second, we fix the sample size and gradually increase number of grid points. We find that accuracy of predictions on noisy data is strongly correlated with number of grid points, however this dependence is weaker for clean data.
\end{enumerate}

\begin{table}
\centering
\begin{tabular}{| c | c | c | c | c |}\toprule \hline
1D inverse transform & \multicolumn{2}{c|}{Noiseless} &\multicolumn{2}{c|}{ Noisy ($13\%$ relative noise)}\\ \hline
Sample Size, $\rm N_x=100$  & $a_0$        & $a_1$        & $a_0$        & $a_1$\\  \hline
$\rm N_{data}=250$  & $0.9953634$  & $0.9977753$  & $0.9609166$  & $0.9570264$\\
$\rm N_{data}=500$  & $0.9979478$  & $0.9988417$  & $0.9644154$  & $0.9640230$\\
$\rm N_{data}=1000$ & $0.9990417$  & $0.9992921$  & $0.9600430$  & $0.9586783$\\
$\rm N_{data}=2000$ & $0.9995410$  & $0.9997107$  & $0.9599427$  & $0.9652383$\\
$\rm N_{data}=4000$ & $0.9994279$  & $0.9994974$  & $0.9603496$  & $0.9661519$\\
$\rm N_{data}=8000$ & $0.9998054$  & $0.9998115$  & $0.9614795$  & $0.9619859$\\ \hline
Grid Points, $\rm N_{data}=1000$  & $a_0$        & $a_1$        & $a_0$        & $a_1$\\  \hline
$\rm N_{x}=20$  & $0.9900532$  & $0.9987560$  & $0.8623348$  & $0.8822680$\\
$\rm N_{x}=40$  & $0.9932568$  & $0.9975166$  & $0.9161125$  & $0.9081806$\\
$\rm N_{x}=80$  & $0.9986574$  & $0.9993274$  & $0.9509870$  & $0.9511483$\\
$\rm N_{x}=160$ & $0.9991550$  & $0.9990234$  & $0.9747287$  & $0.9762977$\\
$\rm N_{x}=320$ & $0.9985649$  & $0.9987451$  & $0.9861375$  & $0.9860783$\\
$\rm N_{x}=640$ & $0.9991842$  & $0.9991606$  & $0.9920950$  & $0.9922520$\\
\hline\bottomrule
\end{tabular}
\caption{$\rm R^2 $ coefficient for predicted Zernike coefficients of the one dimensional Poisson problem at different training sample size and number of grid points.}
\label{tab::tab::R2coeff1D}
\end{table}

\textbf{Two dimensional inverse transform}

We consider an example of variable diffusion coefficients parameterized as $D(x,y)=4 + 2a_2 y + \sqrt{3}a_3(2x^2 + 2y^2 - 1)$, with unknown coefficients randomly chosen in range $a_2,a_3\in [0.25,0.75]$. The equations are solved on a square domain $\Omega \in [-\frac{1}{\sqrt{2}}, \frac{1}{\sqrt{2}}]^2$ governed by the Poisson equation:
\begin{align*}
\nabla\cdot \big([4 + 2a_2 y + \sqrt{3}a_3(2x^2 + 2y^2 - 1)]\nabla u\big) + x + y=0,  &\qquad (x,y)\in \Omega, \\
u_{BC}=\cos(\pi x) \cos(\pi y),  &\qquad (x,y)\in\partial\Omega.
\end{align*}

The encoder is composed of two convolutional layers with $32$ and $64$ channels followed by a $2\times 2$ average pooling layer and a dense layer with $128$ neurons, at the bottleneck there are $2$ neurons representing the two unknowns. All activation functions are \texttt{ReLU} except at the bottleneck that has \texttt{Sigmoid} functions. An \texttt{Adam} optimizer is used on a mean squared error loss function.

\begin{figure}
\centering
\subfigure[Regression quality is $\rm R^2=0.9915683$.]{ \includegraphics[width=0.45\linewidth]{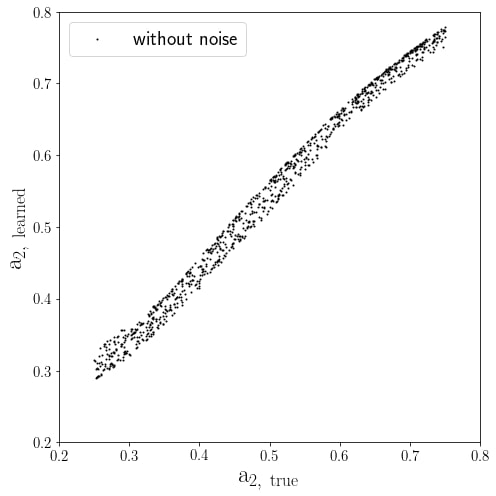} } \quad\quad
\subfigure[Regression quality is $\rm R^2=0.9986852$.]{ \includegraphics[width=0.45\linewidth]{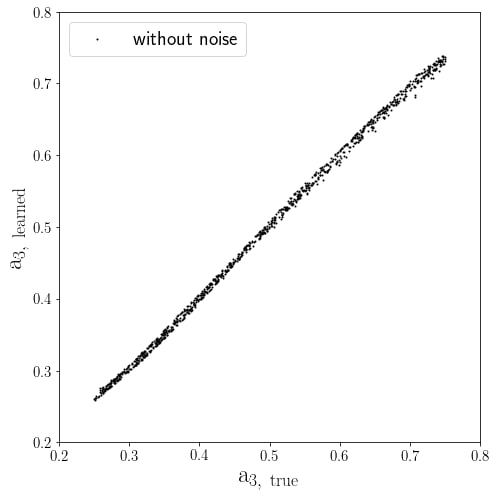} }\quad\quad
\subfigure[Regression quality is $\rm R^2=0.9896654$.]{ \includegraphics[width=0.45\linewidth]{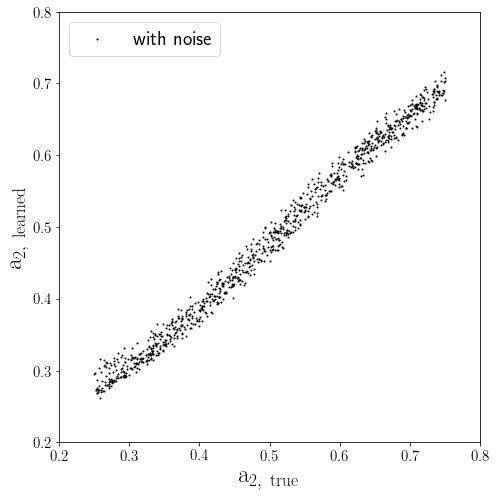} } \quad\quad
\subfigure[Regression quality is $\rm R^2=0.9915149$.]{ \includegraphics[width=0.45\linewidth]{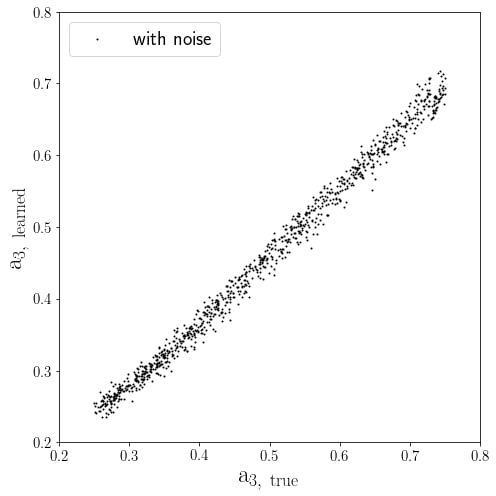} }\quad\quad
\caption{Top row shows performance of BiPDE over $1000$ randomly chosen clean 2D images after 1000 epochs, and the bottom panel shows performance of the same network on noisy images given a zero-mean Gaussian noise with standard deviation $0.025$. Network has $1,852,000$ trainable parameters and training takes $\sim 11$ seconds on a Tesla T4 GPU available on a free Google Colaboratory account.}
\label{fig:2d_inverse}
\end{figure}

We trained BiPDE over $900$ generated solutions with randomly chosen parameters $a_2,a_3$ and validated its performance on 100 independent solution fields for $300$ epochs, evolution of loss function is shown in figure \ref{subfig::lossinv2D}. Then we tested the trained model over another set of $1000$ images with randomly generated diffusion maps independent of the training dataset. Furthermore, we repeated this exercise over $1000$ images with added zero-mean Gaussian noise with standard deviation $0.025$. In each case, the learned parameters are in good agreement with the true values, as illustrated in figure \ref{fig:2d_inverse}. Moreover, we performed a sensitivity analysis on the accuracy of the predictions with respect to sample size. We measured quality of fit by the $\rm R^2$ statistical coefficient. Results are tabulated in table \ref{tab::tab::Two_Unknowns_inverseTran2D} and indicate training over more samples leads to more accurate predictions when applied to clean data, while noisy data do not show a strong improvement by increasing sample size. 

\begin{table}
\centering
\begin{tabular}{| c | c | c | c | c |}\toprule \hline
2D inverse transform & \multicolumn{2}{c|}{$\rm Noiseless$} &\multicolumn{2}{c|}{ Noisy ($13\%$ relative noise)}\\ \hline
Sample Size  & $a_2$        & $a_3$        & $a_2$        & $a_3$\\  \hline
$\rm N_{data}=250$  & $0.9897018$  & $0.9958963$  & $0.9872887$  & $0.9892064$\\
$\rm N_{data}=500$  & $0.9917211$  & $0.9977917$  & $0.9910183$  & $0.9900091$\\
$\rm N_{data}=1000$ & $0.9915683$  & $0.9986852$  & $0.9896654$  & $0.9915149$\\
$\rm N_{data}=2000$ & $0.9940470$  & $0.9993891$  & $0.9909640$  & $0.9883151$\\
$\rm N_{data}=4000$ & $0.9938268$  & $0.9997119$  & $0.9919061$  & $0.9898697$\\
\hline\bottomrule
\end{tabular}
\caption{$\rm R^2 $ coefficient for predicted Zernike coefficients of the two dimensional Poisson problem by increasing training sample size. Number of grid points are fixed at $\rm 30\times 30$.}
\label{tab::tab::Two_Unknowns_inverseTran2D}
\end{table}

\subsection{Dynamic Burger's problem}\label{sec::meshless}
In this section, we demonstrate the applicability of BiPDEs on time-dependent nonlinear partial differential equations, and we use those results to illustrate the consistency and accuracy of the proposed framework. Similar to previous works \cite{Raissi2017PhysicsID}, we consider the nonlinear Burgers' equation in one spatial dimension,
\begin{align}
&\frac{\partial u}{\partial t} + u\frac{\partial u}{\partial x} = \nu \frac{\partial^2 u}{\partial x^2} &x\in [-1,1],~~ t\in [0,1) \label{eq::Burgers}\\
&u(-1,t)=u(1,t)=0   &u(x,0)=-\sin(\pi x)
\end{align}
where $\nu=1/Re$ with $Re$ being the Reynolds number. Notably, Burgers' equation is of great practical significance for understanding evolution equations as it is nonlinear. Burgers' equation has been used as a model equation for the Navier-Stokes equation and by itself can be used to describe shallow water waves \cite{debnath2011nonlinear}, turbulence \cite{burgers1948mathematical}, traffic flow \cite{nagatani2000density}, and many more. 

\begin{itemize}
\item \textbf{Discretization.}
In our design we adopted the $6^{\rm th}$-order compact finite difference scheme proposed by Sari and Gurarslan (2009) \cite{sari2009sixth} for its simplicity of implementation, its high accuracy and because it leads to a linear system with narrow band and subsequently ensures computational efficiency. This scheme combines a tridiagonal\footnote{Tridiagonal systems of equations can be obtained in $\mathcal{O}(N)$ operations.} sixth-order compact finite difference scheme (CFD6) in space with a low-storage third-order accurate total variation diminishing Runge-Kutta scheme (TVD-RK3) for its time evolution (\cite{shu1988efficient}). In particular, the high-order accuracy associated with this discretization provides highly accurate results on coarse grids. This is an important aspect of BiPDEs as a data-efficient inverse solver, which stems from their capacity to seamlessly blend highly accurate and sophisticated discretization methods with deep neural networks.

The first-order spatial derivatives are given at intermediate points by 
\begin{align}
&\alpha u'_{i-1} + u'_i + \alpha u'_{i+1} = b \frac{u_{i+2} - u_{i-2}}{4h} + a\frac{u_{i+1} - u_{i-1}}{2h}, &i=3,\cdots, N-2
\end{align}
where
\begin{align}
&a=\frac{2}{3}(\alpha+2), &b=\frac{1}{3}(4\alpha -1),
\end{align}
and $h=x_{i+1} - x_i$ is the mesh size, with grid points identified by the index $i=1,2,\cdots, N$. For $\alpha=1/3$ we obtain a sixth order accurate tridiagonal scheme. The boundary points (for non-periodic boundaries) are treated by using the formulas \cite{gaitonde1998high,sari2009sixth},
\begin{align*}
u'_1 + 5 u'_{2} &= \frac{1}{h} \bigg[ -\frac{197}{60}u_1 - \frac{5}{12}u_2 + 5u_3 - \frac{5}{3} u_4 + \frac{5}{12} u_5 -\frac{1}{20}u_6\bigg]\\
\frac{2}{11}u'_1 + u'_2 + \frac{2}{11}u'_3 &= \frac{1}{h}\bigg[ -\frac{20}{33} u_1 -\frac{35}{132}u_2 + \frac{34}{33}u_3 - \frac{7}{33}u_4 + \frac{2}{33}u_5 -\frac{1}{132}u_6\bigg] \\
\frac{2}{11}u'_{N-2} + u'_{N-1} + \frac{2}{11}u'_{N}&=\frac{1}{h}\bigg[ \frac{20}{33}u_N + \frac{35}{132}u_{N-1} - \frac{34}{33}u_{N-2} + \frac{7}{33}u_{N-3} - \frac{2}{33}u_{N-4} + \frac{1}{132}u_{N-5} \bigg]\\
5u'_{N-1} + u'_{N}&=\frac{1}{h}\bigg[ \frac{197}{60}u_N + \frac{5}{12} u_{N-1} - 5 u_{N-2} + \frac{5}{3}u_{N-3} -\frac{5}{12}u_{N-4} + \frac{1}{20} u_{N-5} \bigg]
\end{align*}
 This can be easily cast in the matrix form
 \begin{align}
 BU'=AU
 \end{align}
where $U=[u_1, u_2, \cdots, u_N]^T$ is the vector of solution values at grid points. Furthermore, second order derivatives are computed by applying the first-order derivatives twice\footnote{From implementation point of view this is a very useful feature of this scheme, because $A$ and $B$ are constant matrices that do not change through training it is possible to pre-compute them using \texttt{numpy}'s \cite{numpy} basic data structures, and then simply import the derivative operators into \texttt{TensorFlow}'s custom solver layer using  \texttt{tf.convert\_to\_tensor()} command. },
\begin{align}
BU'' = AU'
\end{align}
Burgers' equation are thus discretized as:
\begin{align}
&\frac{dU}{dt}=\mathcal{L}U, &\mathcal{L}U=\nu ~ U'' - U\tens U',
\end{align}
where $\tens$ is the element-wise multiplication operator and $\mathcal{L}$ is a \textit{nonlinear} operator. We use a low storage TVD-RK3 method to update the solution field from time-step $k$ to $k+1$,
\begin{align}
U^{(1)}&=U_k + \Delta t\mathcal{L}U_k\\
U^{(2)}&=\frac{3}{4}U_k  + \frac{1}{4}U^{(1)} + \frac{1}{4}\Delta t \mathcal{L}U^{(1)}\\
U_{k+1}&=\frac{1}{3} U_k + \frac{2}{3}U^{(2)} + \frac{2}{3}\Delta t\mathcal{L}U^{(2)}
\end{align}
with a CFL coefficient of $1$. Note that this method only requires two storage units per grid point, which is useful for large scale scientific simulations in higher dimensions.

\item \textbf{Training protocol.} For training, we first solve Burgers' equation for $M$ time-steps, then we construct two shifted solution matrices that are separated by a single time-step, \emph{i.e.},
\begin{align}
&\mathcal{U}^{-1}=\begin{bmatrix}
U^1 & U^2 & \cdots & U^{M-1}
\end{bmatrix}   &\mathcal{U}^{+1}=\begin{bmatrix}
U^{2} & U^{3} & \cdots & U^M
\end{bmatrix}
\end{align}
Basically, one step of TVD-RK3 maps a column of $\mathcal{U}^{-1}$ to its corresponding column in $\mathcal{U}^{+1}$ given an accurate prediction for the hidden parameter. Hence, a semantic BiPDE is trained with $\mathcal{U}^{-1}$ and $\mathcal{U}^{+1}$ as its input and output respectively. The unknown diffusion coefficient is discovered by the code at the bottleneck of the architecture. 

\item \textbf{Numerical setup.} To enable direct comparison with PINNs, we also consider a second parameter $\gamma$ in Burger's equation. In these current experiments we train for $2$ unknown parameters $(\nu, \gamma)$ in the Burger's equation given by
\begin{align*}
&u_t + \gamma u u_x - \nu u_{xx}=0, & t\in [0,1], ~ ~~ x\in [-1,1]
\end{align*}
Similar to Raissi \emph{et al.} \cite{Raissi2017PhysicsID} we consider $\nu=0.01/\pi$ and $\gamma=1.0$. For completeness we also recall the loss function used in PINN that encodes Burger's equation as a regularization,
\begin{align*}
MSE=\frac{1}{N}\sum_{i=1}^N \bigg\vert u(t^i_u, x_u^i) - u^i \bigg\vert^2 +\frac{1}{N}\sum_{i=1}^N  \bigg\vert u_t(t_u^i, x_u^i) + \gamma u(t_u^i, x_u^i) u_x(t_u^i, x_u^i) - \nu u_{xx}(t_u^i, x_u^i) \bigg\vert^2
\end{align*}
where $(t^i_u, x^i_u, u^i)$ constitute training data with $N=2000$ observation points in the spatio-temporal domain. In this experiment PINN is composed of $9$ layers with $20$ neurons per hidden layer. It is worth mentioning that we are reporting BiPDE results by considering solutions in a narrower time span $t\in [0,0.2]$.

\item \textbf{Architecture.} Obviously, one can choose a single neuron to represent an unknown parameter $\nu$ or $\gamma$ and in a few iterations an approximate value can be achieved. However, our goal is to train a general purpose encoder that is capable of predicting the unknown value from an input solution pair with arbitrary values of $\nu$ and $\gamma$ and without training on new observations (\emph{cf.} see part \ref{subsec::invD}). Therefore, we consider a BiPDE that is composed of a \texttt{conv1D} layer with $128$ filters and a kernel size of $10$ with \texttt{tanh} activation function, followed by \texttt{AveragePooling1D} with a pool size of $2$ that after flattening is fed to two \texttt{Dense} layers with $20$ and $2$ neurons respectively that apply \texttt{Sigmoid} activation function.  We used the \texttt{Adam} optimizer to minimize the mean absolute error measure for the loss function. 

\item \textbf{Accuracy test.} First, we train for two unknowns in Burger's equation, namely $\nu$ and $\gamma$. We perform a sensitivity analysis for $200$ epochs with different numbers of spatial grid points, as well as different time-steps. In each case, we measure the error between the learned values of $\nu$ and $\gamma$ with their true value $\nu_{\rm true}=0.01/\pi$ and $\gamma_{\rm true}=1.0$. Convergence results of this experiment are tabulated in table \ref{tab::Two_Unknowns_NU} and shown in figure \ref{fig:sensitivity}. We find that increasing the number of grid points (\emph{i.e.} the resolution) improves the accuracy up to almost $700$ grid points before accuracy in $\nu$ (but not $\gamma$) starts to deteriorate. Note the decrease in time-step size does not have a significant effect on accuracy unless large number of grid points $N_x>160$ are considered where decreasing time-step clearly improves results.
\begin{figure}[!h]
\centering
\subfigure[Error in $\nu$ - BiPDE with finite difference method.]{ \includegraphics[width=0.45\linewidth]{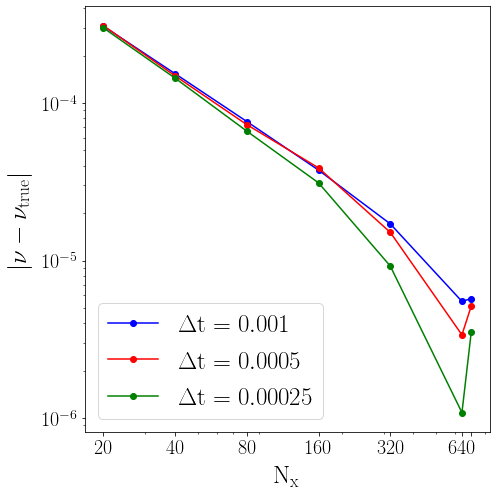} }\quad\quad
\subfigure[Error in $\gamma$ - BiPDE with finite difference method.]{ \includegraphics[width=0.45\linewidth]{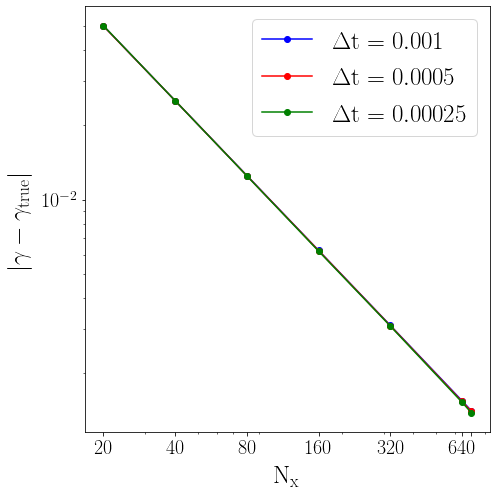} }\quad\quad
\subfigure[Error in $\nu$ - PINN.]{ \includegraphics[width=0.45\linewidth]{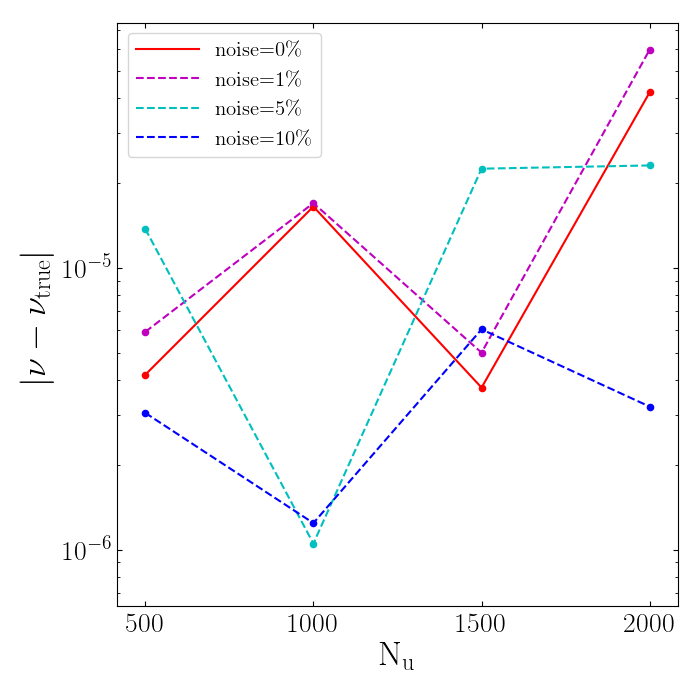}\label{subfig:c} }\quad\quad
\subfigure[Error in $\gamma$ - PINN.]{ \includegraphics[width=0.45\linewidth]{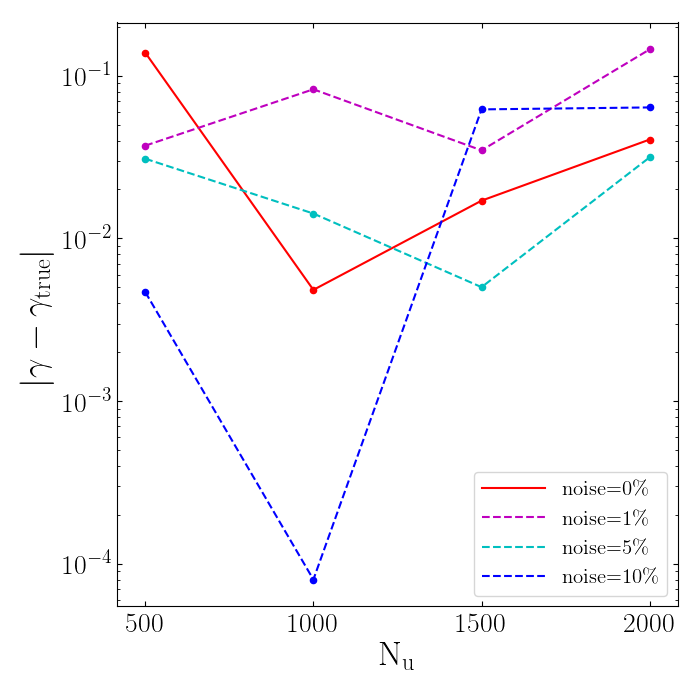} \label{subfig:d}}\quad\quad
\caption{Sensitivity analysis in training both parameters $\gamma$ and $\nu$ with BiPDE (a,b), also results from table 1 of Raissi \etal (2017) \cite{Raissi2017PhysicsID} are shown for comparison (c,d) - note only the solid red line may be compared to BiPDE results where no noise is considered on the solution data. True values are $\nu_{\rm true}=0.01/\pi$ and $\gamma=1.0$. In figure (a) the data points at the right end of $N_x$ axis correspond to $N_x=700$ grid points where the accuracy in the discovered $\nu$ value deteriorates.}
\label{fig:sensitivity}
\end{figure}
\begin{table}
\centering
\begin{tabular}{| c | c | c | c| c | c| c| }\toprule \hline
$\rm \#~ epochs=200$ &  \multicolumn{2}{c}{$\rm \Delta t=0.001$} &  \multicolumn{2}{c}{$\rm \Delta t = 0.0005$} &  \multicolumn{2}{c}{$\rm \Delta t = 0.00025$} \\ \hline
grid size    & $\nu$         & $\gamma$     & $\nu$ & $\gamma$ & $\nu$ & $\gamma$ \\  \hline
$\rm N_x=20$   &  $0.0028751$  &  $0.9500087$ & $0.0028731$ & $0.9500685$  & $0.0028828$    & $0.9499334$ \\
$\rm N_x=40$   &  $0.0030294$  &  $0.9750050$ & $0.0030341$ & $0.9750042$  & $0.0030391$    & $0.9750047$\\
$\rm N_x=80$   &  $0.0031067$  &  $0.9875077$ & $0.0031101$ & $0.9875285$  & $0.0031167$    & $0.9875455$\\
$\rm N_x=160$ &  $0.0031455$  &  $0.9937580$ & $0.0031443$ & $0.9937674$  & $0.0031519$   & $0.9937985$\\
$\rm N_x=320$ &  $0.0031659$  &  $0.9968843$ & $0.0031679$ & $0.9968919$  & $0.0031738$   & $0.9969027$\\
$\rm N_x=640$ &  $0.0031775$  &  $0.9984500$ & $0.0031797$ & $0.9984597$  & $0.0031841$   & $0.9984711$ \\
$\rm N_x=700$ &  $0.0031773$  &  $0.9985866$ & $0.0031779$ & $0.9985945$  & $0.0031865$   & $0.9986123$\\
\hline\bottomrule
\end{tabular}
\caption{Discovering two unknown values of $\nu$ and $\gamma$ in Burger's equation. These values are plotted in figure \ref{fig:sensitivity}.}
\label{tab::Two_Unknowns_NU}
\end{table}

For comparison purposes we report numerical results from table 1 of Raissi \etal (2017) \cite{Raissi2017PhysicsID} in our figures \ref{subfig:c}--\ref{subfig:d}. Here we only presented noise-less results of BiPDE, therefore only the $0\%$ added noise case of PINN is comparable, \emph{i.e.} the solid red line in figures \ref{subfig:c}--\ref{subfig:d}. Even though the two test cases have significant differences and much care should be taken to directly compare the two methods, however BiPDEs have a clear advantage by exhibiting \emph{convergence} in the unknown values, \emph{i.e.} more data means better results.

\begin{table}
\centering
\begin{tabular}{| c | c | c | c|}\toprule \hline
$\rm   <\hat{\nu}>$ & $\rm \Delta t=0.001$ & $\rm \Delta t = 0.0005$ & $\rm \Delta t = 0.00025$ \\ \hline\hline
$\rm N_x=20$   &  $0.0064739$  & $0.0065189$  & $0.0065514$ \\
$\rm N_x=40$   &  $0.0048452$  & $0.0048200$  & $0.0047086$ \\
$\rm N_x=80$   &  $0.0040260$  & $0.0040324$  & $0.0039963$ \\
$\rm N_x=160$  &  $0.0036042$  & $0.0036011$  & $0.0036310$ \\
$\rm N_x=320$  &  $0.0033958$  & $0.0034144$  & $0.0033827$ \\
$\rm N_x=640$  &  $0.0032919$  & $0.0032895$  & $0.0032916$ \\
$\rm N_x=700$  &  $0.0032829$  & $0.0032816$  & $0.0032906$ \\
\hline\bottomrule
\end{tabular}
\caption{Discovering one unknown parameter in Burger's equation, values for $\nu$. }
\label{tab::One_Unknowns_NU}
\end{table}
In a second experiment, we fix the value of $\gamma=1.0$ and only train for the unknown diffusion coefficient $\nu$. Similar to previous test we trained the network for $200$ epochs and figure \ref{fig:sensitivity2} shows the error in the discovered value of $\nu$ at different levels of resolution. In this case decreasing time-step size does not seem to have a significant effect on accuracy. A curious observation is the absolute value of error for $\nu$ is two orders of magnitude more precise when the network is trained for both parameters $\nu$ and $\gamma$ than when only tuning for $\nu$. Again, convergence in unknown parameter is retained in this experiment.
\begin{figure}[ht]
\centering
\includegraphics[width=0.6\linewidth]{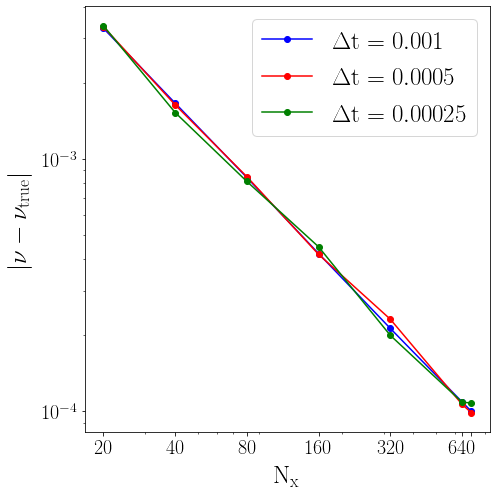} 
\caption{Sensitivity analysis in training only one parameter $\nu$. True value of $\nu_{\rm true}=0.01/\pi$ is sought in Burgers' equation at different levels of resolution. Rightmost data points correspond to $N_x=700$ grid points.}
\label{fig:sensitivity2}
\end{figure}

\end{itemize}

\section{Mesh-less BiPDE: Multi-Quadratic Radial Basis Functions}\label{sec::meshfree}
Not only are direct computations of partial derivatives from noisy data extremely challenging, in many real world applications, measurements can only be made on scattered point clouds. Tikhonov regularization type approaches have been devised to avoid difficulties arising from high sensitivity of differencing operations on noisy data \cite{cullum1971numerical,chartrand2011numerical,stickel2010data}; for neural network based approaches, see \cite{maas2012recurrent,shen2017denoising}. Recently, Trask \etal \cite{trask2019gmls} have proposed an efficient framework for learning from unstructured data that is based on the Generalized Moving Least Squares (GMLS) technique. They show performance of GMLS-Nets to identify differential operators and to regress quantities of interest from unstructured data fields. Another interesting approach had been put forth in the late $80$s by \cite{broomhead1988radial,poggio1990networks} that designed neural networks based on Radial Basis Functions (RBF) to perform functional interpolation tasks. In these networks, the activation function is defined as the radial basis function of interest and the training aims to discover the weights of this network, which interestingly coincide with the coefficients in the corresponding radial basis expansion.

Since the early 70s, RBFs have been used for highly accurate interpolation from scattered data. In particular, Hardy \cite{hardy1971multiquadric} introduced a special kind of RBF called the \textit{multiquadratic} series expansion, that provides superior performance in terms of accuracy, stability, efficiency, simplicity and memory usage \cite{franke1982scattered}. Kansa (1990) \cite{kansa1990multiquadricsI,kansa1990multiquadricsII} pioneered the use of radial basis functions to solve time dependent partial differential equations through deriving a modified multi-quadratic scheme. In 1998, Hon and Mao \cite{hon1998efficient} applied multiquadratics as a spatial discretization method for the nonlinear Burgers' equation and solved it for a wide range of Reynolds numbers (from 0.1 to 10,000). Their scheme was later enhanced to second-order accuracy in time by Xie and Li (2013) \cite{xie2013meshless} via introducing a compact second-order accurate time discretization. Interestingly, the accuracy of these mesh-free methods can be improved by fine-tuning distributions of collocation points or their \textit{shape parameters}. For example, Hon and Mao devised an adaptive point to chase the peak of shock waves, which improved their results. Fortunately, such fine-tuning of parameters can be automated using BiPDE networks; we demonstrate this in this section.

\begin{itemize}
\item \textbf{Discretization.} We chose to blend the second-order accurate method of Xie and Li, briefly described next and we leave further details to their original paper \cite{xie2013meshless}. Initially, one can represent a distribution $u(\mathbf{x})$ in terms of a linear combination of radial basis functions,
\begin{align}
&u(\mathbf{x})\approx \sum_{j=0}^{N_s} \lambda_j \phi_j(\mathbf{x}) + \psi(\mathbf{x}), &\mathbf{x}\in\Omega \subset \mathcal{R}^{dim},\label{eq::RBF}
\end{align}
where $\phi(\mathbf{x})$ is the radial basis function that we adopt,
\begin{align}
&\phi_j(\mathbf{x})=\sqrt{r_j^2 + c_j^2}, &r_j^2 = \vert\vert \mathbf{x} - \mathbf{x}_j\vert\vert_2^2,
\end{align}
and $c_j$ is the \textit{shape parameter} that has been experimentally shown to follow $c_j = Mj + b$ with $j=0, 1, \cdots, N_s$ ($N_s$ is number of seed points). Moreover $M$ and $b$ are tuning parameters. In equation \eqref{eq::RBF}, $\psi(\mathbf{x})$ is a polynomial to ensure solvability of the resulting system when $\phi_j$ is only conditionally positive definite. To solve PDEs, one only needs to represent the solution field with an appropriate form of equation \eqref{eq::RBF}. In the case of Burgers' equation the solution at any time-step $n$ can be represented by
\begin{align}
&u^n(x)\approx \sum_{j=0}^{N_s} \lambda_j^n \phi_j(x) + \lambda_{N_s+1}^n x + \lambda_{N_s+2}^n \label{eq::RBFXL}
\end{align}
over a set of reference points for the basis functions that are given by $x_j=j/N_s$, $j=0,1,\cdots, N_s$. Xie and Li derived the following compact second-order accurate system of equations
\begin{align}
\big[1 + \frac{\Delta t}{2} u_x^n(\hat{x}_j) \big] u^{n+1}(\hat{x}_j) + \frac{\Delta t}{2} u^n(\hat{x}_j) u^{n+1}_x(\hat{x}_j) - \frac{\nu \Delta t}{2}u_{xx}^{n+1}(\hat{x}_j) = u^n(\hat{x}_j) + \frac{\nu \Delta t}{2}u_{xx}^n(\hat{x}_j) \label{eq::XieLi}
\end{align}
over a set of $N_d+1$ distinct collocation points $\hat{x}_j=(1+j)/(N_d+2)$ with $j=0,1,\cdots, N_d$. Two more equations are obtained by considering the left and right boundary conditions $u^{n+1}(x_{L}) = u^{n+1}(x_{R})=0$. Note that spatial derivatives are directly computed by applying derivative operator over equation \eqref{eq::RBFXL}. At every time-step, one solves for the $N+3$ coefficients $\lambda^n_0,\cdots, \lambda^n_{N+2}$, while the spatial components of the equations remain intact (as long as points are not moving). The solution is obtained over the initial conditions given by $u^0(\hat{x}_j)$.

For implementation purposes, we represent the system of equations \eqref{eq::XieLi} in a matrix notation that is suitable for tensorial operations in \texttt{TensorFlow}. To this end, we first write equation \eqref{eq::RBFXL} as
\begin{align}
U^n(\hat{x})= A(\hat{x}) \Lambda^n, \label{eq::linear}
\end{align}
where 
\begin{align}
&U^n_{(N_d+1)\times 1}  = \begin{bmatrix}
u^n(\hat{x}_0) \\ u^n(\hat{x}_1) \\ \vdots \\ u^n(\hat{x}_{N_d})
\end{bmatrix},    & \Lambda^n_{(N_s+1)\times 1} = \begin{bmatrix}
\lambda_0^n \\ \lambda_1^n \\ \vdots \\ \lambda_{N_s}^n
\end{bmatrix},
\end{align}
\begin{align}
&\bigg[A_{ij}(\hat{\mathbf{x}})\bigg]_{(N_d+1) \times (N_s+1)}=\bigg[\phi_{j}(\hat{x}_i) - \phi_{j}(x_L) - \frac{ \phi_{j}(x_R) - \phi_j(x_L)}{x_R -x_L}(\hat{x}_i - x_L)  \bigg],
\end{align}
with $i=0,1,\cdots, N_d$ and $j=0,1,\cdots, N_s$. Note that we already injected the homogeneous boundary conditions into equation \eqref{eq::linear}. Therefore, equation \eqref{eq::XieLi} can be written as,
\begin{align}
\bigg[ A + (\mathbf{g}_x~\mathbf{1}^T)\tens A + (\mathbf{g}~ \mathbf{1}^T)\tens A_x - \frac{\nu \Delta t}{2}A_{xx}\bigg] \Lambda^{n+1} = \bigg[A + \frac{\nu \Delta t}{2}A_{xx}\bigg] \Lambda^n, \label{eq::RBFequat}
\end{align}
where $\mathbf{1}^T=[1,~ 1, ~\cdots, ~1]_{1\times(N_s+1)}$, $\tens$ is component-wise multiplication, and
\begin{align}
&\mathbf{g}= \frac{\Delta t}{2}A \Lambda^n,             &\mathbf{g}_x= \frac{\Delta t}{2}A_x \Lambda^n, \\
&(A_x)_{ij}=\phi'_{j}(\hat{x}_i) - \frac{ \phi_{j}(x_R) - \phi_j(x_L)}{x_R -x_L},            &(A_{xx})_{ij}=\phi''_{j}(\hat{x}_i).
\end{align}

Note that in case of training for two parameters $(\nu,\gamma)$, expression for $\mathbf{g}$ in equation \ref{eq::RBFequat} needs to be modified by letting $\mathbf{g}=\frac{\gamma \Delta t}{2}A\Lambda^n$.

\item \textbf{Architecture.} Note that both the collocation points and the interpolation seed points can be any random set of points within the domain and not necessarily a uniform set of points as we chose above. In fact, \textit{during training we allow BiPDE to find a suitable set of interpolation points as well as the shape parameters on its own}. The input data is calculated using aforementioned finite difference method over uniform grids and later interpolated on a random point cloud to produce another sample of solutions on unstructured grids for training. Thus, in our architecture the last layer of the encoder has $2N_s + 1$ neurons with \texttt{sigmoid} activation functions representing the $2N_s$ shape parameters and seed points, as well as another neuron for the unknown diffusion coefficient. Note that for points to the left of origin, in the range $x\in [-1, 0]$, we simply multiplied the output of $N_s$ activation functions by $``-1''$ within the solver layer (because output of \texttt{Sigmoid} function is always positive). We use the mean squared error between data and predicted solution at time-step $n+p$ as the loss function. We used the \texttt{Adam} optimizer to minimize the loss function.

\item \textbf{Training protocol.} As in the previous case, we apply successive steps of MQ-RBF scheme to march the input data forward to a future time-step. Not surprisingly, we observed that taking higher number of steps improves the results because erroneous guess of the diffusion coefficient leads to more pronounced discrepancy after longer periods of time. Hence, we map the data $\mathcal{U}^{-p}$ to $\mathcal{U}^{+p}$, which is $p$ time-steps shifted in time,
\begin{align}
&\mathcal{U}^{-p}=\begin{bmatrix}
U^1, U^2, \cdots, U^{M-p}
\end{bmatrix}     &\mathcal{U}^{+p}=\begin{bmatrix}
U^{1+p}, U^{2+p}, \cdots, U^{M}
\end{bmatrix} 
\end{align}
In our experiments a value of $ p=10$ was sufficient to get satisfactory results at the absence of noise. However, at the presence of Gaussian noise and for smaller values of the diffusion coefficient (such as for $\nu_{\rm true}=0.01/\pi$) we had to increase the shift parameter to $ p=100$.

\item \textbf{Numerical setup.} Once again, we let $\nu_{\rm true}=0.01/\pi\approx 0.00318$ and integrate Burgers' equation up to $t_f=0.2$ with a fixed time-step of $\Delta t=0.001$. We use the finite difference method of the previous section to generate the datasets. We then interpolate the solution on 80 data points, uniformly distributed in the range $(-1,1)$ with 20 interpolation seed points. For this experiment, we set the batch size to $1$. We trained the network using \texttt{Adam} optimizer. The results after 50 epochs are given in figure \ref{fig:results_RBF_proc}.
\begin{figure}
\centering
\subfigure[True solution generated by finite differences (input data).]{ \includegraphics[height=0.37\linewidth]{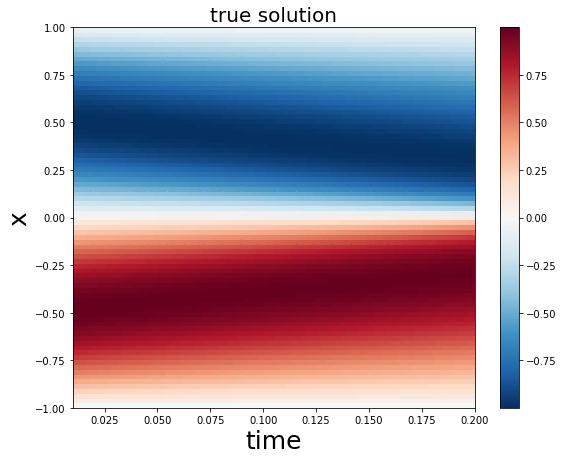} } \quad\quad
\subfigure[Learned solution generated by MQ-RBF BiPDE (output data).]{ \includegraphics[height=0.37\linewidth]{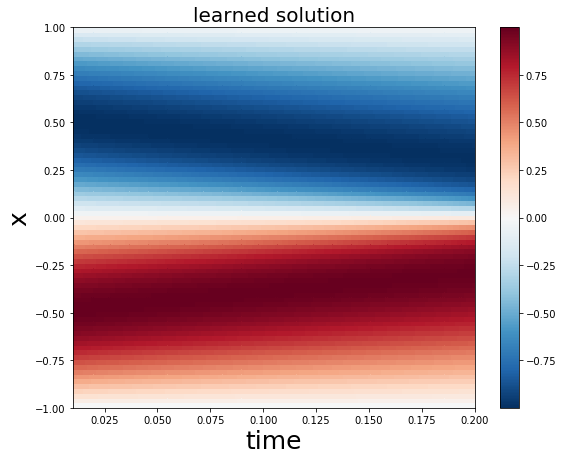} }\quad\quad
\subfigure[Error in solution.]{ \includegraphics[height=0.4\linewidth]{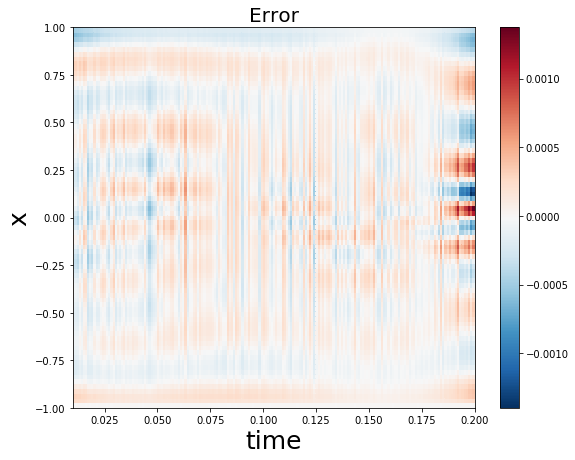} }\quad\quad
\subfigure[Discovered seeds and shape parameters.]{ \includegraphics[height=0.4\linewidth]{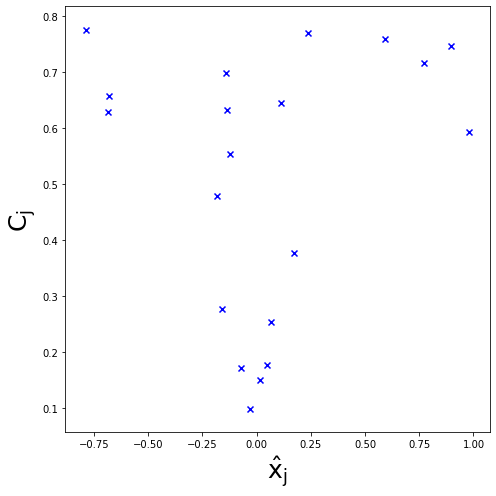} }\quad\quad
\subfigure[Distribution of diffusion coefficients.]{ \includegraphics[width=0.45\linewidth]{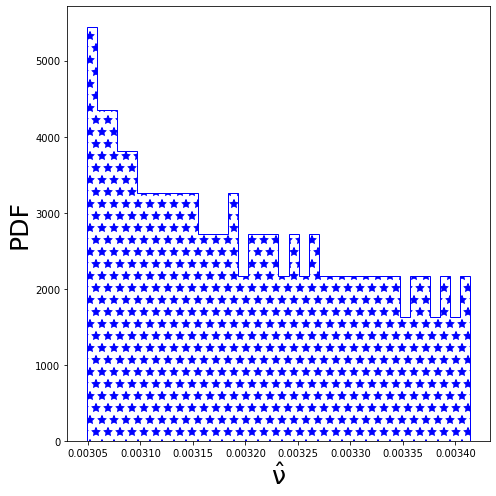} }\quad\quad
\subfigure[Evolution of mean squared error during training.]{ \includegraphics[height=0.45\linewidth]{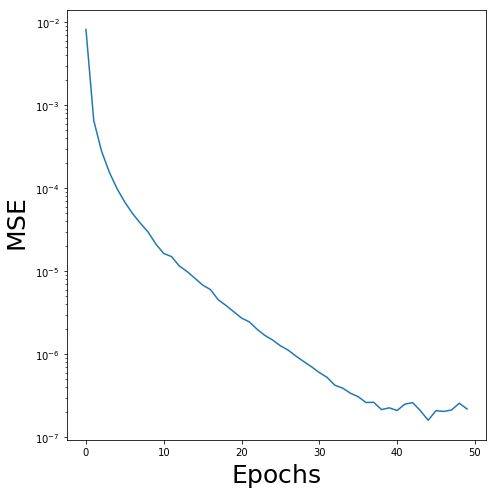} }\quad\quad
\caption{Results of applying the RBF-BiPDE to Burgers' equation with a true diffusion coefficient of $\nu_{\rm true}=0.003183$. The average value of the predicted diffusion coefficients is $\hat{\nu}=0.00320$.}
\label{fig:results_RBF_proc}
\end{figure}

Interestingly, for every pair of input-output, the network discovers a distinct value for the diffusion coefficient that provides a measure of uncertainty for the unknown value. We report the average value of all diffusion coefficients as well as the probability density of these values. We observe that for all pairs of solutions, the predicted value for the diffusion coefficient is distributed in the range $0.00305\le \hat{\nu} \le 0.00340$ with an average value of $<\hat{\nu}>=0.00320$, which is in great agreement with the true value, indeed with $0.6\%$ relative error. Interestingly, we observe that the BiPDE network has learned to concentrate its interpolation seed points around the origin where the solution field varies more rapidly. Furthermore, around $x=\pm 0.5$, the interpolation points are more sparse, which is in agreement with the smooth behavior of the solution field at these coordinates. Therefore, this network may be used as an automatic shock tracing method to improve numerical solutions of hyperbolic problems with shocks and discontinuities as was shown by Hon and Mao.

\item \textbf{Resilience to noise on unstructured grids.} We consider several cases to assess robustness to noise. In each case, we pick $80$ \textit{randomly} distributed points along the $x$-axis and linearly interpolate the solution field on this set of points. Then, we add a Gaussian noise with a given standard deviation. This noisy and unstructured data field is then fed into the MQ-RBF based BiPDE of this section. We use a batch size of 10, with $10\%$ of each sample for validation during training. A summary of our results follows:
\begin{enumerate}
\item Let $\nu_{\rm true}=0.1/\pi$, $p=10$, $N_d=80$, $N_s=20$, $\Delta t=0.001$, and consider a Gaussian noise with a standard deviation of $1\%$. After $100$ epochs, we obtain the results in figure \ref{fig:results_RBF_proc2}.
\begin{figure}
\centering
\subfigure[True solution generated by finite differences and with added noise. Solution is interpolated on a random grid.]{ \includegraphics[height=0.35\linewidth]{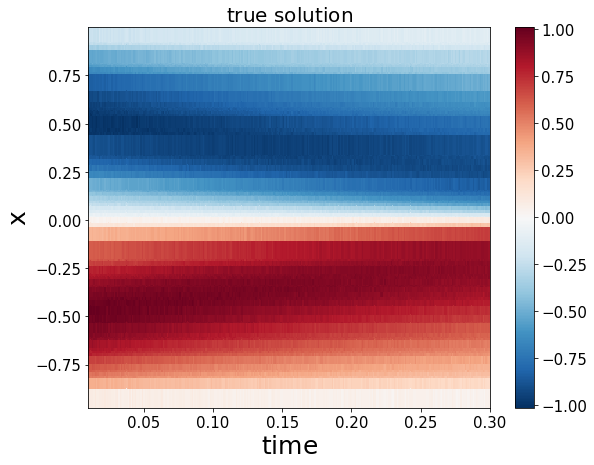} } \quad\quad
\subfigure[Learned solution generated by MQ-RBF BiPDE (output data).]{ \includegraphics[height=0.35\linewidth]{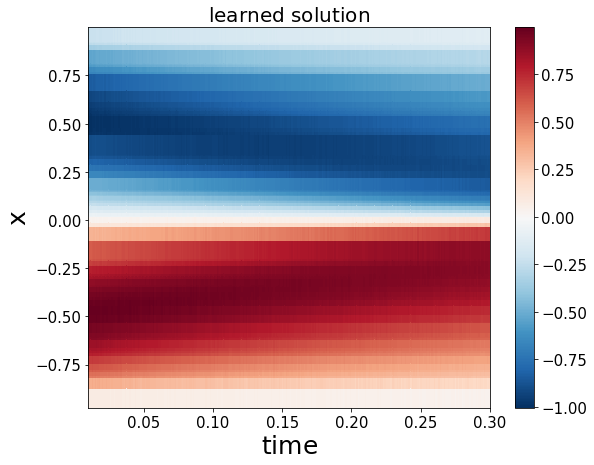} }\quad\quad
\subfigure[Error in solution.]{ \includegraphics[height=0.4\linewidth]{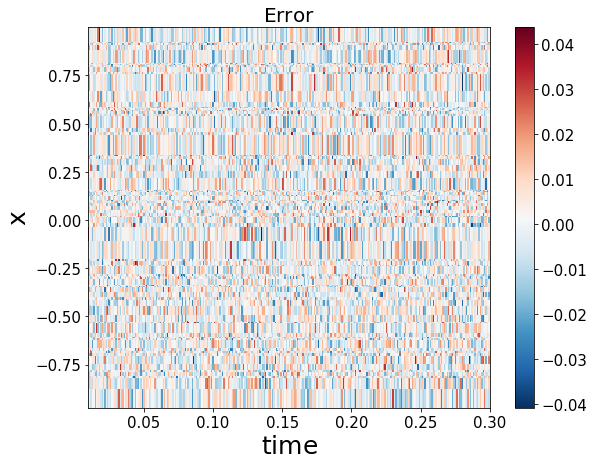} }\quad\quad
\subfigure[Discovered seeds and shape parameters. Error bars indicate one standard deviation.]{ \includegraphics[height=0.4\linewidth]{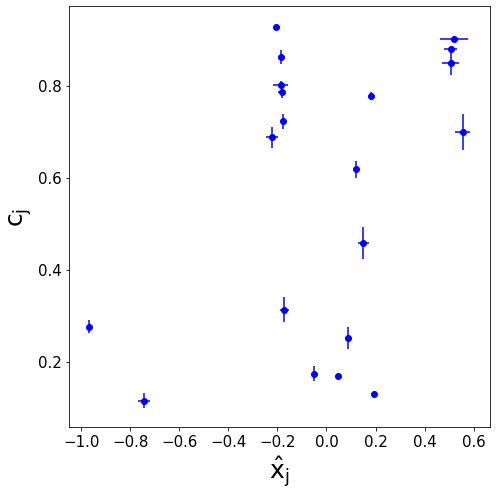} }\quad\quad
\subfigure[Probability density of diffusion coefficients.]{ \includegraphics[width=0.45\linewidth]{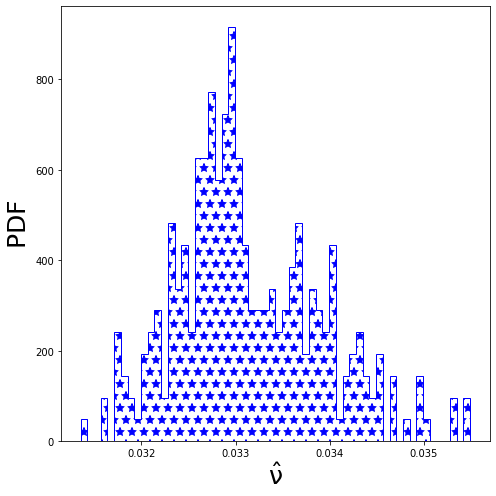} }\quad\quad
\subfigure[Evolution of mean squared error versus number of epochs.]{ \includegraphics[width=0.45\linewidth]{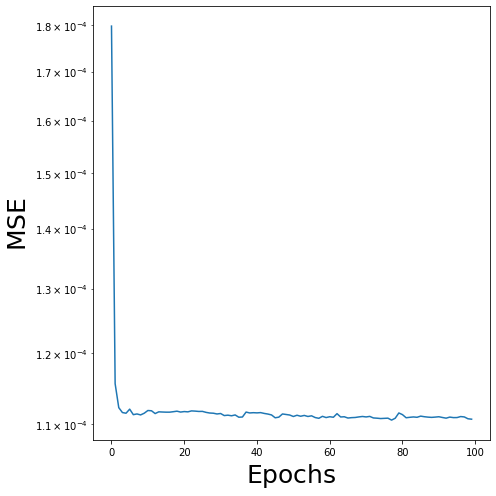} }\quad\quad
\caption{Results of applying the RBF-BiPDE to Burgers' equation with a true diffusion coefficient of $\nu_{\rm true}=0.03183$. The average value of the predicted diffusion coefficients is $\hat{\nu}=0.0331$. The data is provided on a scattered point cloud with added Gaussian noise with $1\%$ standard deviation.}
\label{fig:results_RBF_proc2}
\end{figure}

\item Let $\nu_{\rm true}=0.1/\pi$, $p=100$, $N_d=200$, $N_s=20$, $\Delta t=0.001$, and consider a Gaussian noise with a standard deviation of $5\%$. After $150$ epochs, we obtain the results in figure \ref{fig:results_RBF_proc3}.
\begin{figure}
\centering
\subfigure[True solution generated by finite differences and with added noise. Solution is interpolated on a random grid.]{ \includegraphics[height=0.35\linewidth]{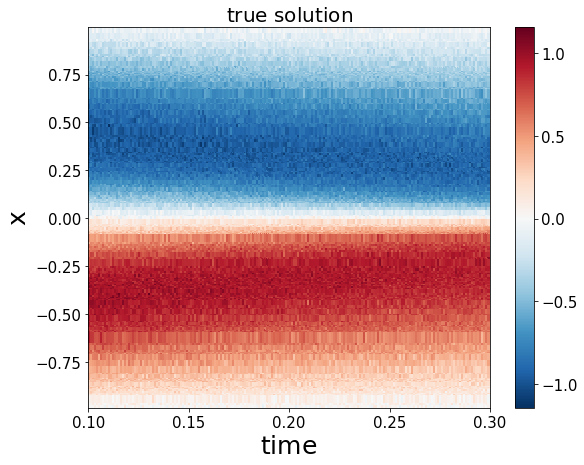} } \quad\quad
\subfigure[Learned solution generated by MQ-RBF BiPDE (output data).]{ \includegraphics[height=0.35\linewidth]{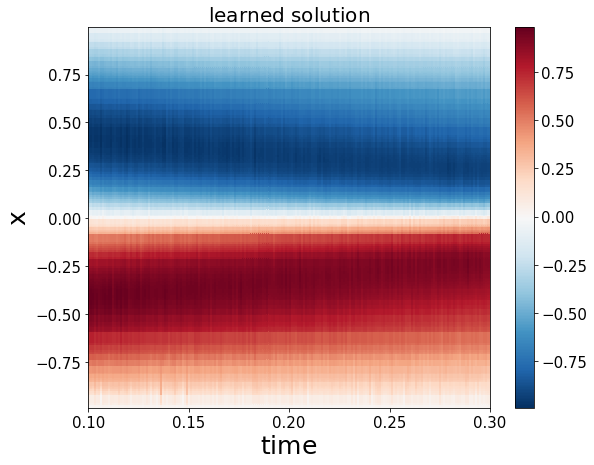} }\quad\quad
\subfigure[Error in solution.]{ \includegraphics[height=0.4\linewidth]{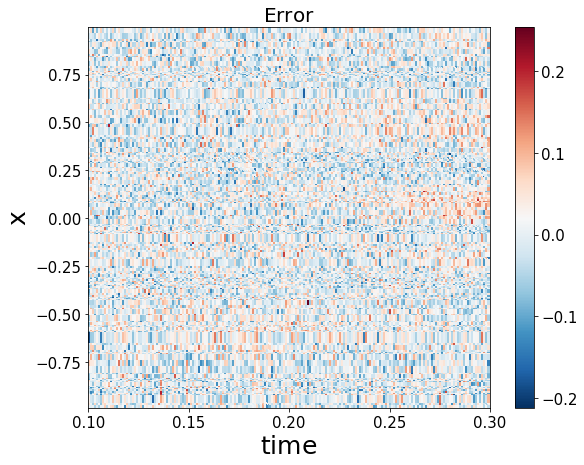} }\quad\quad
\subfigure[Discovered seeds and shape parameters. Error bars indicate one standard deviation.]{ \includegraphics[height=0.4\linewidth]{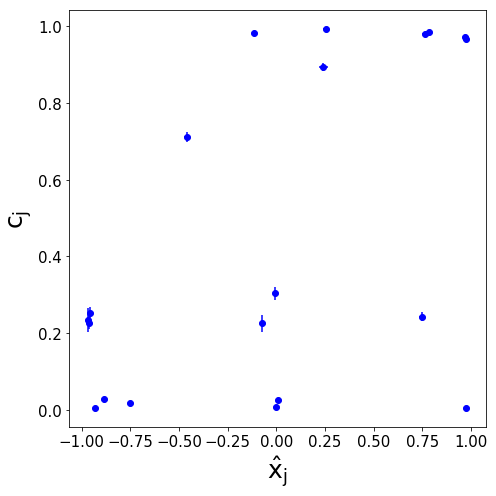} }\quad\quad
\subfigure[Probability density of diffusion coefficients.]{ \includegraphics[width=0.45\linewidth]{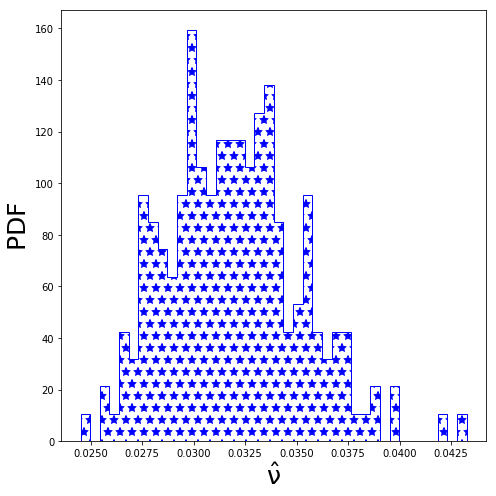} }\quad\quad
\subfigure[Evolution of mean squared error versus number of epochs.]{ \includegraphics[width=0.45\linewidth]{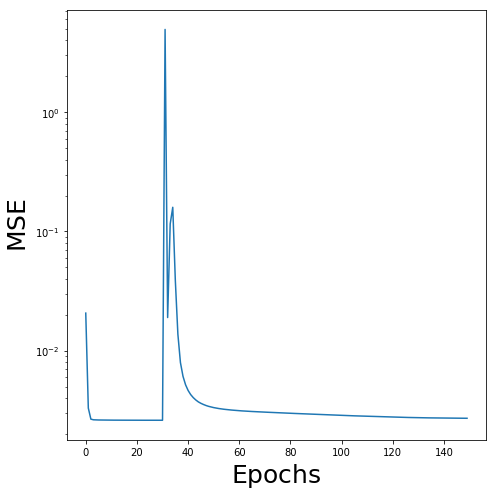} }\quad\quad
\caption{Results of applying the RBF-BiPDE to Burgers' equation with a true diffusion coefficient of $\nu_{\rm true}=0.03183$. The average value of the predicted diffusion coefficients is $\hat{\nu}=0.03160$. The data is provided on a scattered point cloud with added Gaussian noise with $5\%$ standard deviation.}
\label{fig:results_RBF_proc3}
\end{figure}

\item Let $\nu_{\rm true}=0.01/\pi$, $p=100$, $N_d=80$, $N_s=20$, $\Delta t=0.001$, and consider a Gaussian noise with a standard deviation of $1\%$. After $200$ epochs, we obtain the results in figure \ref{fig:results_RBF_proc4}.
\begin{figure}
\centering
\subfigure[True solution generated by finite differences and with added noise. Solution is interpolated on a random grid.]{ \includegraphics[height=0.35\linewidth]{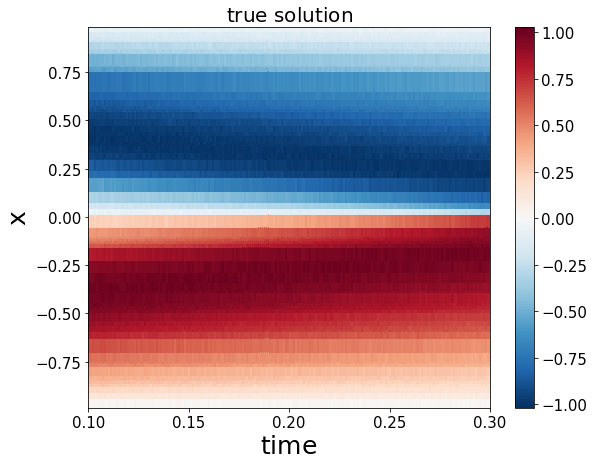} } \quad\quad
\subfigure[Learned solution generated by MQ-RBF BiPDE (output data).]{ \includegraphics[height=0.35\linewidth]{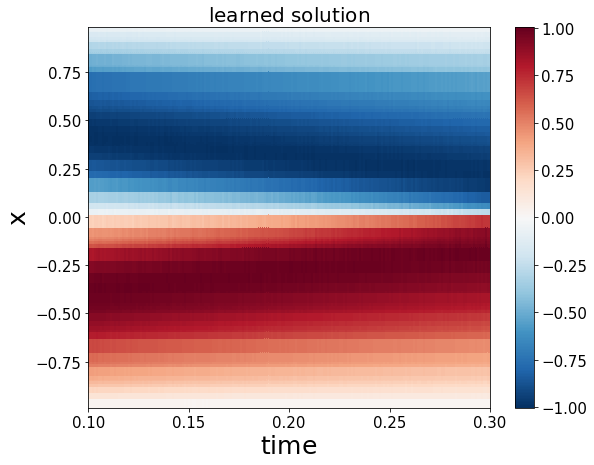} }\quad\quad
\subfigure[Error in solution.]{ \includegraphics[height=0.4\linewidth]{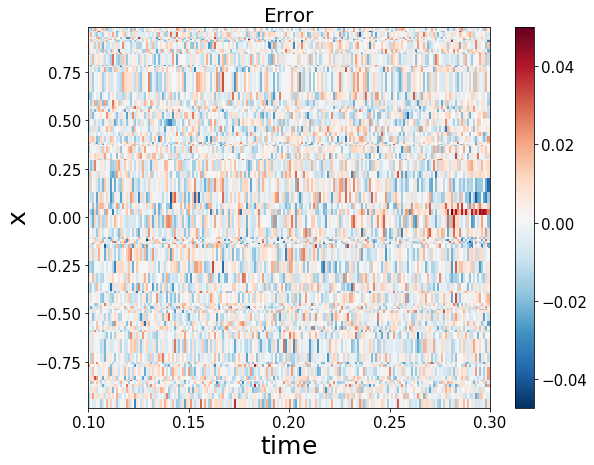} }\quad\quad
\subfigure[Discovered seeds and shape parameters. Error bars indicate one standard deviation.]{ \includegraphics[height=0.4\linewidth]{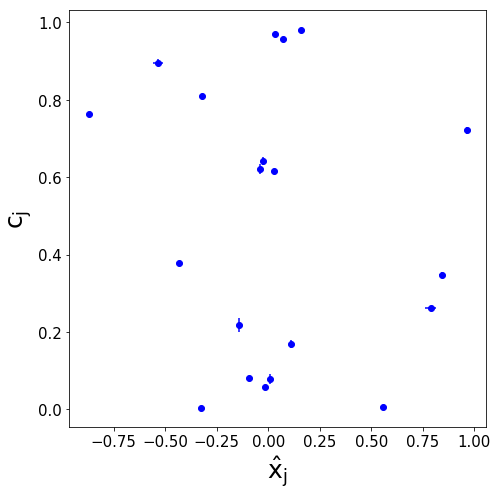} }\quad\quad
\subfigure[Probability density of diffusion coefficients.]{ \includegraphics[width=0.45\linewidth]{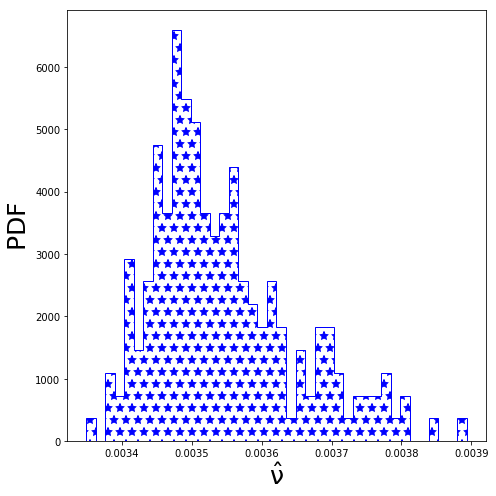} }\quad\quad
\subfigure[Probability density of diffusion coefficients.]{ \includegraphics[width=0.45\linewidth]{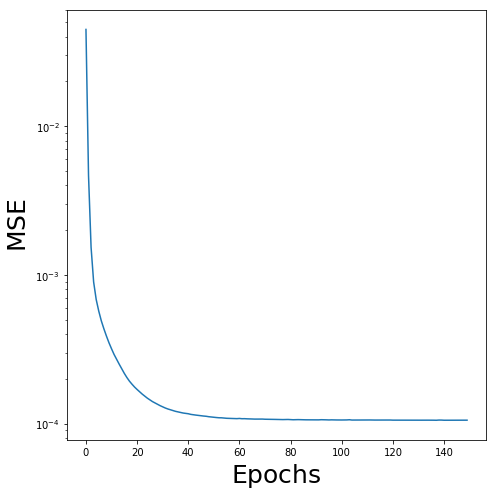} }\quad\quad
\caption{Results of applying the RBF-BiPDE to Burgers' equation with a true diffusion coefficient of $\nu_{\rm true}=0.003183$. The average value of the predicted diffusion coefficients is $\hat{\nu}=0.00352$. The data is provided on a scattered point cloud with added Gaussian noise with $1\%$ standard deviation.}
\label{fig:results_RBF_proc4}
\end{figure}

\item Let $\nu_{\rm true}=0.01/\pi$, $p=100$, $N_d=200$, $N_s=20$,  $\Delta t=0.001$, and consider a Gaussian noise with a standard deviation of $5\%$. After $150$ epochs, we obtain the results in figure \ref{fig:results_RBF_proc5}.
\begin{figure}
\centering
\subfigure[True solution generated by finite differences and with added noise. Solution is interpolated on a random grid.]{ \includegraphics[height=0.35\linewidth]{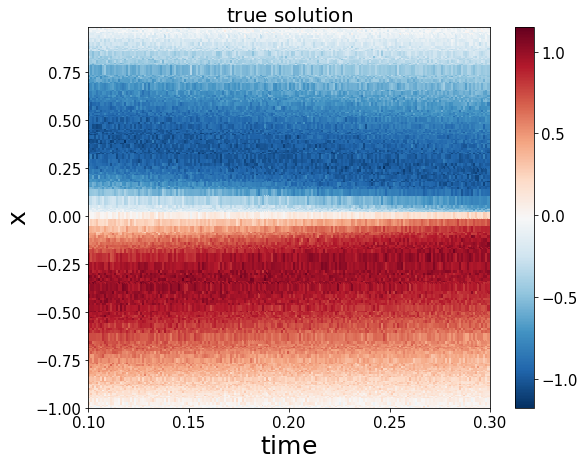} } \quad\quad
\subfigure[Learned solution generated by MQ-RBF BiPDE (output data).]{ \includegraphics[height=0.35\linewidth]{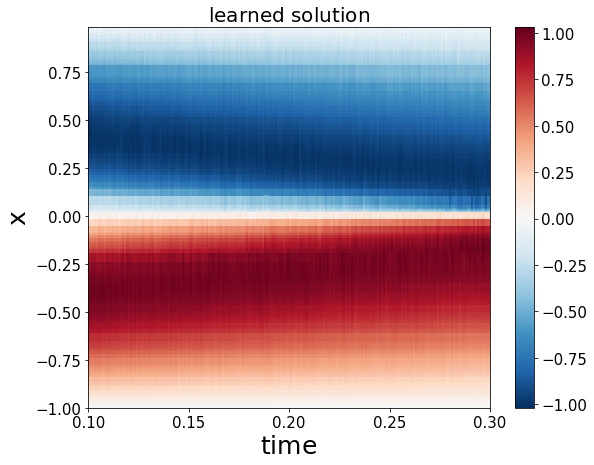} }\quad\quad
\subfigure[Error in solution.]{ \includegraphics[height=0.4\linewidth]{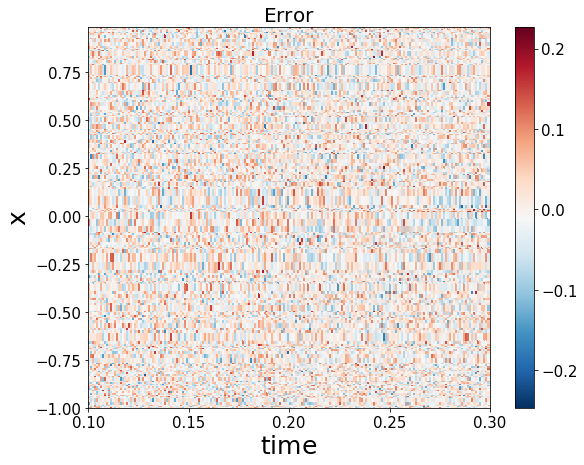} }\quad\quad
\subfigure[Discovered seeds and shape parameters. Error bars indicate one standard deviation.]{ \includegraphics[height=0.4\linewidth]{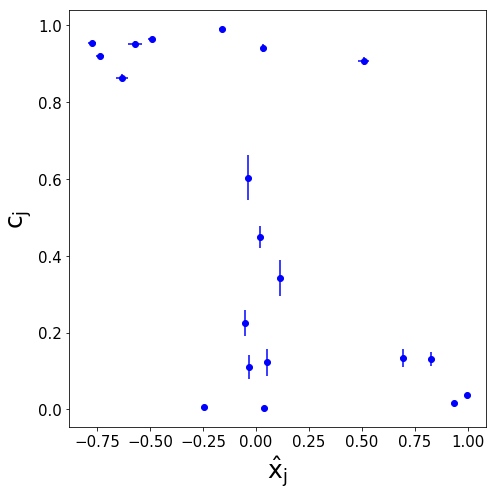} }\quad\quad
\subfigure[Probability density of diffusion coefficients.]{ \includegraphics[width=0.45\linewidth]{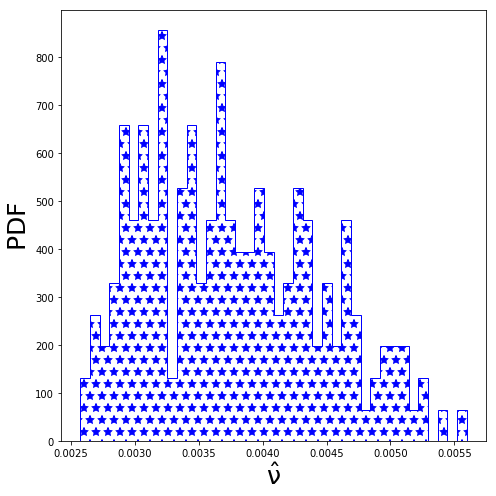} }\quad\quad
\subfigure[Probability density of diffusion coefficients.]{ \includegraphics[width=0.45\linewidth]{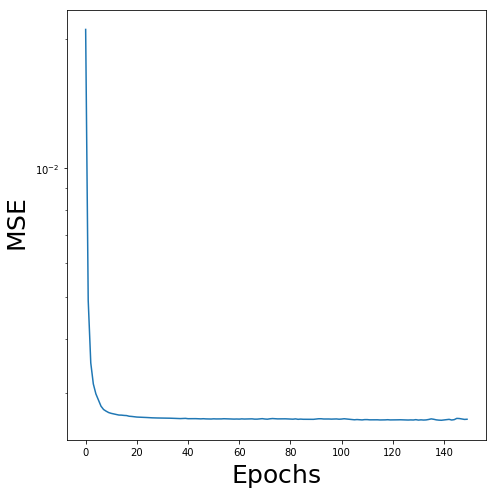} }\quad\quad
\caption{Results of applying the RBF-BiPDE to Burgers' equation with a true diffusion coefficient of $\nu_{\rm true}=0.003183$. The average value of the predicted diffusion coefficients is $\hat{\nu}=0.003677$. The data is provided on a scattered point cloud with added Gaussian noise with $5\%$ standard deviation.}
\label{fig:results_RBF_proc5}
\end{figure}
\end{enumerate}
We observe that this architecture is generally robust to noise. However, at higher noise values require more tuning of hyperparameters, as well as longer training.

\item \textbf{Accuracy tests.}  We report the values of the discovered diffusion coefficients in the Burgers' equation for different grid sizes and different time-steps. We use the same setting as that detailed in the numerical setup part in this section. Particularly, the interpolation seeds are determined by the network and the training data is on a uniformly distributed set of points computed by the finite difference method of the previous section. We consider three different time-steps, $\Delta t=0.001, ~0.0005, ~0.00025$, and two diffusion coefficients of $\nu_{\rm true}=0.01/\pi,~0.1/\pi$ over grids of size $N_x=80,~ 160$. At each time-step, for all experiments with different grid sizes, we stop the training when the mean squared error in the solution field converges to a constant value and does not improve by more epochs; this roughly corresponds to $50, ~ 25,~12$ epochs for each of the training time-steps, respectively. This indicates that by choosing smaller time steps less number of epochs are needed to obtain the same level of accuracy in the unknown parameter. Furthermore, we use an \texttt{Adam} optimizer with a learning rate of 0.001.

The results of the accuracy tests are tabulated in tables \ref{tab::nu1}--\ref{tab::nu2}. We observe, for all experiments, that the discovered coefficient is in great agreement with the true values. Due to adaptivity of the interpolation seed points and their shape parameters for different experiments, the observed error values do not seem to follow the trend of traditional finite difference methods, as depicted in previous sections. This could also be due to lower order of accuracy of the MQ-RBF method, \textit{i.e.} being a second-order accurate method, compared to the higher-order accurate finite difference method used in the previous section.

\begin{table}
\centering
\begin{tabular}{| c | c | c | c|}\toprule \hline
$  <\hat{\nu}>$ & $\Delta t=0.001$ & $\Delta t = 0.0005$ & $\Delta t = 0.00025$ \\ \hline
$\#~ epochs$    & $50$         & $25$     & $12$\\  \hline
$N_x=80$   &  $0.03173\pm 3.4\times 10^{-4}$  & $0.03196 \pm 4.2\times 10^{-4}$  & $0.03188 \pm 2.8\times 10^{-4}$ \\
$N_x=160$ &  $0.03186 \pm 5.8\times 10^{-5}$  & $0.03191\pm 3.6\times 10^{-4}$  & $0.03137\pm 1.2\times 10^{-4}$ \\ \hline\bottomrule
\end{tabular}
\caption{Discovered values of the diffusion coefficient for $\nu_{\rm true}=0.03183$ at different time-steps and grid sizes. }
\label{tab::nu1}
\end{table}

\begin{table}
\centering
\begin{tabular}{| c | c | c | c|}\toprule \hline
$  <\hat{\nu}>$ & $\Delta t=0.001$ & $\Delta t = 0.0005$ & $\Delta t = 0.00025$ \\ \hline
$\#~ epochs$    & $50$         & $25$     & $12$\\  \hline
$N_x=80$   &  $0.003326\pm 5.1\times 10^{-5}$  & $0.003162\pm 2.2\times 10^{-4}$  & $0.003155\pm 1.2 \times 10^{-4}$ \\
$N_x=160$ &  $0.003264\pm 1.0\times 10^{-4}$  & $0.003151\pm 1.3\times 10^{-4}$  & $0.003192\pm 1.2\times 10^{-4}$ \\ \hline \bottomrule
\end{tabular}
\caption{Discovered values of the diffusion coefficient for $\nu_{\rm true}=0.003183$ obtained with different time-steps and grid sizes.}
\label{tab::nu2}
\end{table}

\end{itemize}

\subsection{Learning the inverse transform}\label{subsec::invD}
\begin{figure}
\centering
\includegraphics[width=1.1\linewidth]{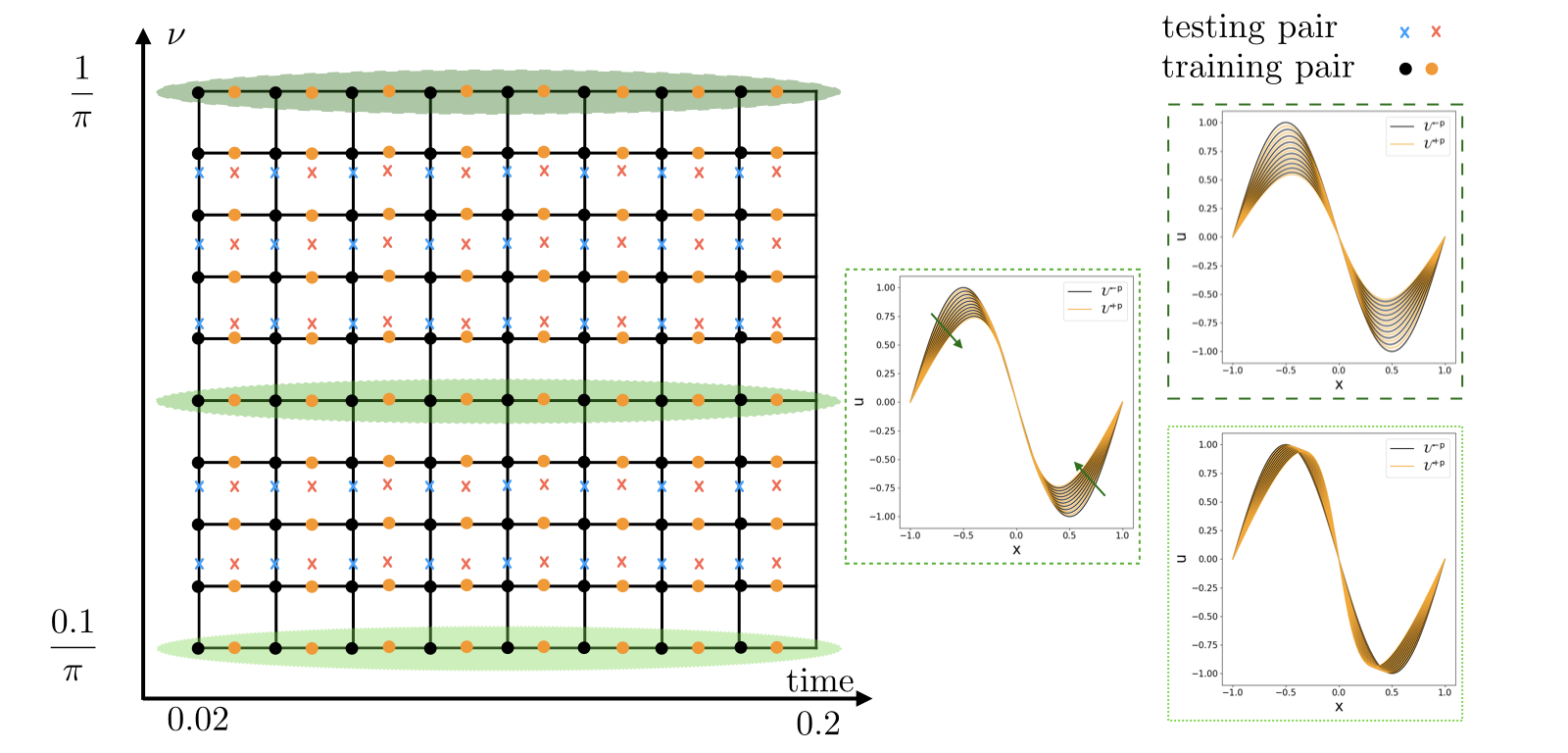}
\caption{Topology of data points for training and testing of the semantic BiPDE. Along the $\nu$ dimension, we depict $10$ (out of $50$) of the selected data points, while along the time dimension we illustrate the actual $8$ data points. Training pairs of $\mathcal{U}^{-p}$ and  $\mathcal{U}^{+p}$ are color coded by black and orange dots, respectively; testing pairs are depicted by blue and red crosses. On the right panel, we illustrate the training data for three nominal values of the diffusion coefficient, highlighted by green shades. Green arrows indicate the direction of time.}
\label{fig::semanticBiPDE_dyn}
\end{figure}

As we emphasized before, a feature of BiPDE is to produce self-supervised pre-trained encoder models for inverse differential problems that are applicable in numerous applications where hidden values should be estimated in real-time. We train an encoder over a range of values $\nu\in [0.1/\pi,~ 1/\pi]$ and assess the performance of the trained model on new data with arbitrarily chosen $\nu$ values. We choose 50 diffusion coefficients that are distributed uniformly in this range, then integrate the corresponding Burgers' equation up to $t_f=0.2$ with a constant time-step of $\Delta t=0.0005$ on a grid with $N_x=100$ grid points using the aforementioned finite difference method. There are $4000$ time-steps in each of the $50$ different realizations of Burgers' equation. For a fixed value of $p=20$, we draw  $10$ solution pairs for each value of $\nu$ at uniformly distributed time instances and discard the first two instances to improve convergence of the network. Hence, the training data uniformly samples the space of solutions over a $8\times 50$ grid of $(t, \nu)$, as illustrated in figure \ref{fig::semanticBiPDE_dyn}. We use the resulting $400$ pairs in training of a semantic BiPDE, with $320$ pairs used for training and $80$ pairs for validation. 

\textbf{Architecture.} Given an arbitrary input, the signature of the hidden physical parameters will be imprinted on the data in terms of complex patterns spread in space and time. We use a CNN layer as a front end unit to transform the input pixels to internal image representations. The CNN unit has $32$ filters with kernel size of $5$. The CNN is followed by max pooling with pool size of $2$, which is then stacked with another CNN layer of $16$ filters and kernel size of $5$ along with another max pooling layer. The CNN block is stacked with two dense layers with $100$ and $41$ neurons, respectively. CNN and dense layers have \texttt{ReLU} and \texttt{Sigmoid} activation functions, respectively. Overall, there are $42,209$ trainable parameters in the network. Conceptually, the CNN extracts features on every snapshot that characterizes the evolution of the solution field through time-steps with a proper physical parameter. This parameter is enforced to be the diffusion coefficient through the PDE solver decoder stage. We train this network for $500$ epochs using an \texttt{Adam} optimizer.

\textbf{Resilience to noise.} Even though the encoder is trained on ideal datasets, we demonstrate a semantic BiPDE still provides accurate results on noisy datasets. In contrast to other methods, we pre-train the network in a self-supervised fashion on clean data and later we apply the trained encoder on unseen noisy data\footnote{Note that the network could also be trained on noisy data as we showed before; however training would take longer in that case.}.
\begin{figure}
\centering
\subfigure[Performance of encoder on training data set.]{ \includegraphics[height=0.35\linewidth]{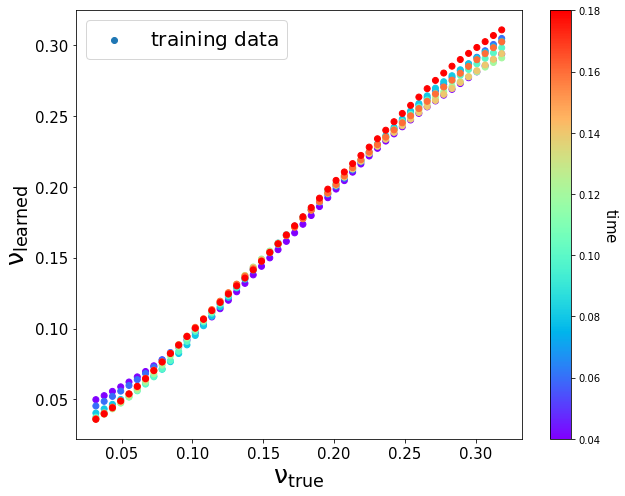} } \quad\quad
\subfigure[Distribution of interpolation points and shape parameters discovered by the network.]{ \includegraphics[height=0.35\linewidth]{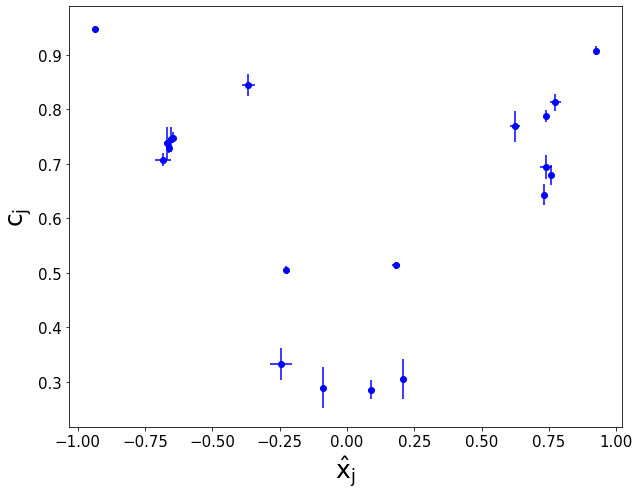} } \quad\quad
\subfigure[Performance of the encoder on unseen data.]{ \includegraphics[height=0.35\linewidth]{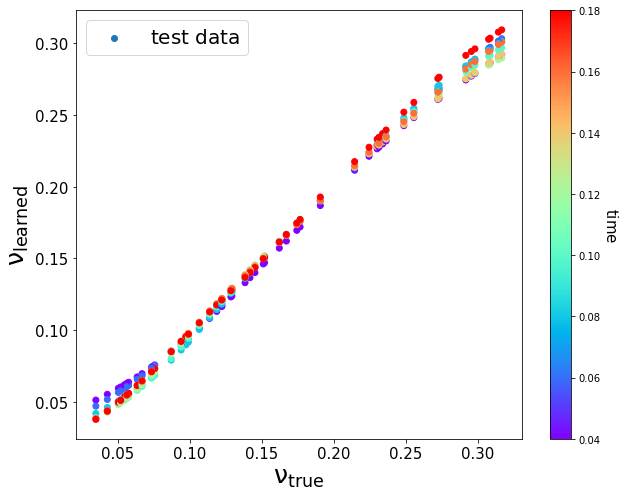} } \quad\quad
\subfigure[Performance of the encoder on unseen data with Gaussian noise with standard deviation $0.01$.]{ \includegraphics[height=0.35\linewidth]{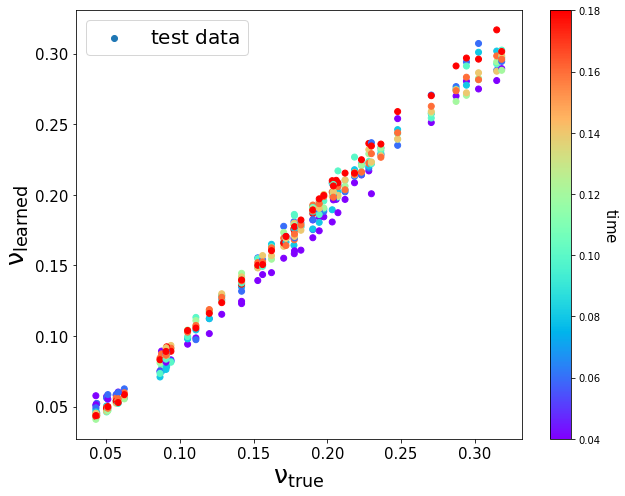} } \quad\quad
\caption{Semantic autoencoder learns how to discover hidden variables from pairs of solutions. These results are obtained after $500$ epochs on $50$ data points along the $\nu$-axis.}
\label{fig:semanticBiPDE_dyn}
\end{figure}
In figure \ref{fig:semanticBiPDE_dyn}, we provide the performance of this network on training as well as on unseen clean/noisy data-sets. Furthermore, the network determines optimal parameters of the MQ-RBF method by evaluating interpolation seed points as well as their corresponding shape parameters to obtain the best approximation over \textit{all} input data.

\section{Conclusion}
We introduced BiPDE networks, a natural architecture to infer hidden parameters in partial differential equations given a limited number of observations. We showed that this approach is versatile as it can be easily applied to arbitrary static or nonlinear time-dependent inverse-PDE problems. We showed the performance of this design on multiple inverse Poisson problems in one and two spatial dimensions as well as on the non-linear time-dependent Burgers' equation in one spatial dimension. Moreover, our results indicate BiPDEs are robust to noise and can be adapted for data collected on unstructured grids by resorting to traditional mesh-free numerical methods for solving partial differential equations. We also showed the applicability of this framework to the discovery of inverse transforms for different inverse-PDE problems. 

There are many areas of research that could be further investigated, such as considering diffusion maps with discontinuities across subdomains, using more sophisticated neural network architectures for more complex problems, using higher-order numerical solvers and finally tackle more complicated governing PDE problems with a larger number of unknown fields or in higher dimensions.

\section*{Acknowledgment}
This research was supported by ARO W911NF-16-1-0136 and ONR N00014-17-1-2676.

\newpage
\section*{References}
\bibliographystyle{abbrv}
\addcontentsline{toc}{section}{\refname}
\bibliography{references}

\end{document}